\documentclass[12pt]{article} \usepackage{amsfonts,amssymb,latexsym} 
\newtheorem{theorem}{Theorem}     

\newtheorem{lemma}[theorem]{Lemma} \newtheorem{corollary}[theorem]{Corollary} 
 
\def\beq{\begin{equation}}  \def\eeq{\end{equation}}
\def\bb{\begin{eqnarray*}}  \def\ee{\end{eqnarray*}}
\def\b{\begin{eqnarray}}    \def\e{\end{eqnarray}}
\def\1{\hbox{\rm\setbox1=\hbox{1}\copy1\kern-.5\wd1 I}}
\def\D{{\hbox{\rm\setbox1=\hbox{I}\copy1\kern-.45\wd1 D}}}
\def\E{{\hbox{\rm\setbox1=\hbox{I}\copy1\kern-.45\wd1 E}}}
\def\N{\hbox{\rm\setbox1=\hbox{I}\copy1\kern-.55\wd1 N}}
       \def\L{{\cal L}}
  
\def\p{\hbox{\rm\setbox1=\hbox{I}\copy1\kern-.45\wd1 P}}
\def\R{\hbox{\rm\setbox1=\hbox{I}\copy1\kern-.45\wd1 R}} 
\def\<{{^{_<}}}   \def\>{{^{_>}}}   \def\st{\stackrel{d}{=}}
\def\nin{\in{\!\!\!\!\!/}\,}        \def\bib{\vspace{-2mm}\bibitem}

\def\tp{\rightarrow_{\!\!\!\!\!\!_p}\ \,}  \def\pr{\noindent{\bf Proof}\ } 
  \def\t{\bar {t\,}}
\def\z{\zeta}   \def\l{\lambda} \def\vp{\varphi}   \def\ve{\varepsilon} 
\def\a{\alpha}    \def\t{\tau} 
\def\Z{{\bf Z}_+}   \def\({\left(}  \def\){\right)} \def\d{d_{_{TV}}\!} 

\begin{document} 

 \author{S.Y.Novak\\ {\small MDX University London}\\ {\small\it version 2}} 
 \date{} 
 \title{\vspace*{-2.9cm}Poisson approximation} \maketitle\vspace*{-0.6cm} 
 \begin{abstract} 
   We overview results on the topic of Poisson approximation that are missed in existing surveys. 
	The main attention is paid to the problem of Poisson approximation to the distribution of a sum of Bernoulli and, more generally, non-negative integer-valued random variables. 

   We do not restrict ourselves to a particular method, and overview the whole range of issues including the general limit theorem, estimates of the accuracy of approximation, asymptotic expansions, etc. 
	Related results on the accuracy of compound Poisson approximation are presented as well. 
 
  We indicate a number of open problems and discuss directions of further research. \end{abstract} 
 

	
     { \noindent {\bf Content}\\ 
\noindent\ref{p1}.\,\ \ Weak convergence to a Poisson law\hspace*{\fill}\pageref{p1}\\ 
\ref{p1}.1 Independent random variables\hspace*{\fill}\pageref{p1.1}\\ 
\ref{p1}.2 Dependent Bernoulli random variables\hspace*{\fill}\pageref{p1.2}\\ 
\ref{2}.\,\ \ Accuracy of Poisson approximation\hspace*{\fill}\pageref{2}\\
   2.1 Independent Bernoulli random variables\hspace*{\fill}\pageref{2.1}\\ 
	 2.2 Dependent Bernoulli random variables\hspace*{\fill}\pageref{2.2}\\  
   2.3 Independent integer-valued random variables\hspace*{\fill}\pageref{2.3}\\
	 2.4 Dependent integer-valued random variables\hspace*{\fill}\pageref{2.4}\\ 
	 2.5 Asymptotic expansions\hspace*{\fill}\pageref{2.6}\\ 
	 2.6 Sum of a random number of random variables\hspace*{\fill}\pageref{2.5}\\ 
\ref{3}.\,\ \ Applications\hspace*{\fill} \pageref{3}\\ 
   \ref{3}.1 Long head runs\hspace*{\fill}\pageref{3.1}\\
	 \ref{3}.2 Long match patterns\hspace*{\fill}\pageref{3.2}\\ 
\ref{5}.\,\ \ Compound Poisson approximation\hspace*{\fill}\pageref{5}\\ 
   \ref{5}.1 CP limit theorem\hspace*{\fill}\pageref{5}\\ 
   \ref{5}.2 Accuracy of CP approximation\hspace*{\fill}\pageref{5.1}\\ 
\ref{5}.3 CP approximation to $\,{\bf B}(n,p)$\hspace*{\fill}\pageref{5.2}\\ 
\ref{4}.\,\ \ Poisson process approximation\hspace*{\fill}\pageref{4}\\ 
	 6.\,\ \ References\hspace*{\fill}\pageref{6}\\ } 

			\section{Weak convergence to a Poisson law}\label{p1}  

   Poisson approximation appears natural in situations where one deals with a large number of rare events. The topic has attracted a considerable body of research. It has important applications in insurance, extreme value theory, reliability theory, mathematical biology, etc. (cf. \cite{BalK,BHJ,H67,LLR,N11}). 
	However, existing surveys are surprisingly sketchy, and miss not only a number of results obtained during the last three decades but even some classical results going back to 1930s. 
	
	The paper aims to fill the gap. We present a comprehensive list of results on the topic of Poisson approximation, and formulate a number of open problems. 
	Related results on the topic of compound Poisson approximation are presented as well. 
	
				\subsection{Weak convergence to a Poisson law}\label{p1.1}

  We denote by $\,{\bf\Pi}(\l)\,$ a Poisson law with parameter $\,\l$. 
	
	The following  Poisson limit theorem is due to Gnedenko \cite{Gn38} and Marcinkiewicz \cite{M38}. Hereinafter multiplication is superior to division. 

Let $\,\{X_{n,1},...,X_{n,k_n}\}_{n\ge1},$ where $\,\{k_n\}\,$ is a non-decreasing sequence of natural numbers, be a triangle array of {\it independent} random variables (r.v.s). 

   Random variables $\,\{X_{n,k}\}\,$ are called {\it infinitesimal} if 
 \beq                                                \label{Inf}
 \lim_{n\to\infty}\max_{k\le k_n} \p(|X_{n,k}|\!>\!\ve) \to 0 \qquad(\forall\ve\!>\!0). \eeq 
   Denote $\,B_\ve = (-\ve;\ve) \cup (1\!-\!\ve;1\!+\!\ve),$ 
$$\,S_n = X_{n,1} + ... + X_{n,k_n}. $$

  \begin{theorem} \label{GM} {\rm\cite{Gn38, M38}} 
If $\,\{X_{n,k}\}\,$ are infinitesimal r.v.s, then  
	\beq 										\label{Pois} 
	\L(S_n) \Rightarrow {\bf\Pi}(\l) \qquad (\exists\l\!>\!0) 
	\eeq 
as $\,n\to\infty\,$ if and only if for any $\,\ve\!\in\!(0;1)$, as $\,n\!\to\!\infty,$ 
 \b \label{GM-1} 
 && \sum_k \p(|X_{n,k}\!-\!1|\!<\!\ve) \to \l, 		\label{GM-2}\\ 
 && \sum_k \p(X_{n,k} \!\nin\! B_\ve) \to 0,\ \  
    \sum_k \E X_{n,k} \1\{|X_{n,k}|\!<\!\ve\} \to 0, \label{GM-3}\\ 
 && \sum_k \(\E X_{n,k}^2 \1\{|X_{n,k}|\!<\!\ve\} - 
 \E^2 X_{n,k}\1\{|X_{n,k}|\!<\!\ve\}\) \to 0. 	  \label{GM-4}
 \e \end{theorem} 

  The following corollary presents necessary and sufficient conditions for the weak convergence of a sum of independent and identically distributed (i.i.d.) non-negative integer-valued r.v.s to a Poisson random variable. 

	Let $\,\N\,$ denote the set of natural numbers, and let $\,\Z\!:=\!\N\cup\!\{0\}$. 

  \begin{corollary} \label{CorGM} If $\,\{X_{n,1},...,X_{n,k_n}\}_{n\ge1}\,$ is a triangle array of independent random variables taking values in $\,\Z\,$ such that $\,X_{n,i}\st X_{n,1}\ (1\!\le\!i\!\le\!k_n)$, then (\ref{Pois}) holds if and only if 
 \beq																							\label{P10} 
 k_n\p(X_{1,n}\!=\!1) \to\l\ \ \hbox{\rm and}\ \, 
 \p(X_{n,1}\!\ge\!2)/\p(X_{n,1}\!\ge\!1) \to0.  \eeq \end{corollary} 

	Note that (\ref{P10}) yields 
	$\,\p(X_{n,1}\!\ge\!1) \sim \p(X_{n,1}\!=\!1)\,$ as $\,n\!\to\!\infty$. 
	
	The second relation in (\ref{P10}) means $\,X_{n,1}'\tp1\,$ as $\,n\!\to\!\infty$, where r.v. $\,X_{n,1}'\,$ has the distribution 
	$\,\L(X_{n,1}') = \L(X_{n,1}|X_{n,1}\!\ne\!0).$\\

  In the case of Bernoulli $\,{\bf B}(p_{n,k})\,$ random variables relations (\ref{GM-3}) and (\ref{GM-4}) trivially hold, (\ref{GM-2}) means 
	$$ \sum_k p_{n,k} \to \l \qquad(n\to\infty), \eqno(\ref{GM-2}^*) $$ 
while (\ref{Inf}) states that $\,\max_kp_{n,k}\to0\,$ as $\,n\to\infty$. The latter together with (\ref{GM-2}$^*$) is equivalent to 
  $$ \sum_k p_{n,k}^2 \to0 \qquad(n\!\to\!\infty). \eqno(\ref{Inf}^*) $$ 
Thus, conditions $\,(\ref{Inf}^*)\,$ and $\,(\ref{GM-2}^*)\,$ are necessary and sufficient for the weak convergence (\ref{Pois}).\\ 

  \noindent{\bf Example 1.1}. Let $\,\{X_{n,1},...,X_{n,n}\}\,$ be i.i.d. random variables with the distribution 
	$$ 
	\p(X\!=\!0)=1\!-\!\l/n\!-\!1/n^{1.5},\ \p(X\!=\!1)=\l/n,\ 
	\p(X\!=\!n) = 1/n^{1.5} \qquad(\l\!>\!0). 
	$$
Then (\ref{Inf}) and (\ref{P10}) hold, hence $\,\L(S_n) \Rightarrow {\bf\Pi}(\l)$. Note that $\,\E S_n$ $\!\to{\!\!\!\!\!\!/}$ $\l.$ \hspace*{\fill}$\Box$\\ 

	The proof of Theorem \ref{GM} can be found in \cite{GK54}. 

   A compound Poisson limit theorem (weak convergence of $\,\L(S_n)\,$ to a compound Poisson law, where $\,S_n\,$ is a sum of i.i.d. random variables that are equal to 0 with a large probability) has been given by Khintchin (\cite{Khi}, ch. 2.3). 

		\subsection{Dependent Bernoulli random variables}\label{p1.2} 

   The topic of Poisson approximation to the distribution of a sum of {\it dependent} Bernoulli r.v.s has applications in extreme value theory, reliability theory, etc. (cf. \cite{BalK,BHJ,LLR,N11}).  

   Let $\,\{X_{n,1},...,X_{n,n}\}_{n\ge1}\,$ be a triangle array of 0-$\!$1 random variables such that sequence $\,X_{n,1},...,X_{n,n}\,$ is stationary for each $\,n\!\in\!\N$. 	For instance, in extreme value theory one often has 
	$$ X_{n,k} = \1\{Y_k\!>\!u_n\}, $$ 
where $\,\{Y_i,i\!\ge\!1\}\,$ is a {\it stationary} sequence of random variables and $\,\{u_n\}\,$ is a sequence of ``high'' levels. 
   The special case where $\,\{Y_i,i\!\ge\!1\}\,$ is a moving average is related to the topic of the Erd\"os--R\'enyi partial sums (cf. \cite{N11}, ch. 2). 

   Let $\,{\cal F}_{l,m}(\t)\,$ be the $\,\sigma$--field generated by the events 
$\,\{ X_{n,i}\},$ $\,l\!\le\!i\!\le\!m$. 
   Set 
	\bb  
\alpha_n(l) \!&=&\! \sup|\,\p(AB)\!-\!\p(A)\p(B)|,\ \ \vp(l) = \sup|\p(B|A)-\p(B)|,\\ 
	\beta_n(l)  \!&=&\!  \sup \E\sup_B |\p(B|{\cal F}_{1,m})\!-\!\p(B)|, 
	\ee  
where the supremum is taken over $\,m\!\ge\!1,\,A\!\in\!{\cal F}_{1,m}(\t),$ $\,B\!\in\!{\cal F}_{m+l+1,n}\,$ such that $\,\p(A)\!>\!0$. 
   Conditions involving mixing coefficients $\,\alpha_n(\cdot), \beta_n(\cdot), \vp_n(\cdot)\,$ are slightly weaker than those involving traditional mixing coefficients $\,\alpha(\cdot), \beta(\cdot), \vp(\cdot)$.\\ 

   {\it Condition $\,\Delta\,$} is said to hold if $\,\a_n(l_n)\to0\,$ for some sequence $\,\{l_n\}\,$ of natural numbers such that $\,1\ll l_n\ll n.$\\ 

   {\it Class} $\,{\cal R}$. If $\,\Delta\,$ holds, then there exists a sequence $\,\{r_n\}\,$ of natural numbers such that 
 \beq                                                       \label{rr}
 n\gg r_n\gg l_n \gg1,\ nr_n^{-1} \a_n^{2/3}(l_n) \to0 \qquad(n\!\to\!\infty) 
 \eeq 
(for instance, one can take $\,r_n = \big[\sqrt{n\max\{l_n;n\a_n(l_n)\}}\,\big]$). 
   We denote by $\,{\cal R}\,$ the class of all such sequences $\,\{r_n\}$.\\ 

   Set $$ S_n = X_{n,1}+...+X_{n,n},\ \ \l_n = \E S_n. $$ 
Let $\,\z_{r,n}\,$ be a r.v. with the distribution 
	\beq \label{z} \L(\z_{r,n}) = \L(S_r\!\mid\!S_r\!>\!0). \eeq 
  In extreme value theory $\,\L(\z_{r,n})\,$ is known as the cluster size distribution.

 \begin{theorem} \label{Dep} Assume condition $\,\Delta$. If, as $\,n\!\to\!\infty$, 
 \beq \label{dep} S_n\Rightarrow \pi_\l \qquad(\exists\l\!>\!0),\eeq 
then 
 \beq \label{1} \z_{r,n} \tp1 \qquad(n\!\to\!\infty) \eeq 
for any sequence $\,\{r\!=\!r_n\}\,$ obeying (\ref{rr}). 

 If there exists the limit 
 \beq 																	\label{13.1} 
 \lim_{n\to\infty} \p(X_{n,1}\!=\!...\!=\!X_{n,n}\!=\!0)=e^{-\l} \qquad(\exists\l\!>\!0)  \eeq 
and (\ref{1}) holds for some $\,\{r\!=\!r_n\}\!\in\!{\cal R},$ then $\,S_n\Rightarrow \pi_\l$. \end{theorem}

  Theorem \ref{Dep} generalises Corollary \ref{CorGM} to the case of dependent $\a$-mixing r.v.s. 
	
   Condition (\ref{13.1}) is an analogue of (\ref{GM-2}); it  means that $\,\p(X_{n,i}\!\ne\!0)\,$ are ``properly small''. 
	
   Condition (\ref{1}) prohibits asymptotic clustering of rare events. 
   In the case of independent r.v.s taking values in $\,\Z\,$ assumption (\ref{1}) means $\,X_{n,1}'\tp1\,$ as $\,n\!\to\!\infty$, where r.v. $\,X_{n,1}'\,$ has the distribution $\,\L(X_{n,1}') = \L(X_{n,1}|X_{n,1}\!\ne\!0).$\\

  {\bf Remark 2.1}. The following condition $\,(D')\,$\index{condition $(D')$} has been widely used in extreme value theory (cf. \cite{LLR,N11}): 
 $$
 \lim_{n\to\infty} n\sum_{i=1}^{r-1} \p\(X_{n,i+1}\!\ne\!0, X_{n,1}\!\ne\!0\)=0 
 \eqno{(D')} 
 $$ 
for any sequence $\,\{r\!=\!r_n\}\,$ such that $\,n\!\gg\!r_n\!\gg\!1$. 
   Condition $\,(D')\,$ means that there is no asymptotic clustering of extremes. 
It was introduced by Loynes \cite{Lo}.  
   
	Closely related is the following condition 
 $$
 \lim_{n\to\infty} \sum_{i=1}^{r-1} \p\(X_{n,i+1}\!\ne\!0 | X_{n,1}\!\ne\!0\) = 0. \eqno{(\tilde D')} 
 $$ 
	If conditions $\,\Delta\,$ and (\ref{13.1}) hold, then $\,(D')\,$ is {\it equivalent} to $(\tilde D')$. 
	
	Indeed, one can check that $\,\Delta\,$ and (\ref{13.1}) yield 
 \beq                                                       \label{13.6}
 \p(S_r\!>\!0) \sim \l r/n \qquad(n\!\to\!\infty)
 \eeq	
(cf. (\ref{2b2}) below). Denote $\,p = \p(X_{n,1}\!\ne\!0)$. Then 
  $$ \l r/n \sim \p(S_r\!>\!0) \le rp,\ \l\!+\!o(1) \le np. $$ 
Hence ($D'$) $\Rightarrow$ ($\tilde D'$). 

  By Bonferroni's inequality, 
	\bb
	\l r/n \!&\sim&\! \p(S_r\!>\!0)\ \ge\ r\p(X_{n,1}\!\ne\!0) - 
	\p\Big( \cup_{1\le i<j\le r} \{X_{n,i}\!\ne\!0, X_{n,j}\!\ne\!0\}\Big) \\ 
	\!&\ge&\! rp - rp \sum_{i=1}^{r-1} \p\( X_{n,i+1}\!\ne\!0 | X_{n,1}\!\ne\!0 \). 
	\ee 
Therefore, 
 \b \label{rp} 
 1 \!&\ge&\! \p(S_r\!>\!0)/rp \ge 1-\sum_{i=1}^{r-1}\p\(X_{n,i+1}\!\ne\!0|X_{n,1}
   \!\ne\!0\), \\ \label{bda} \l\!+\!o(1)  \!&\ge&\! \Big( 1 - \sum_{i=1}^{r-1} 
	 \p\(X_{n,i+1}\!\ne\!0 | X_{n,1}\!\ne\!0\) \!\Big) np. 
 \e 
   Thus, $\,np\,$ is bounded away from 0 and above, and $\,(D')\,$ is equivalent to $\,(\tilde D')$.\\

   {\bf Remark 2.2}. Condition (\ref{1}) is weaker than $\,(D')$: if conditions $\,\Delta\,$ and (\ref{13.1}) hold, then $\,(D')\,$ entails (\ref{1}). 
	Indeed, $\,\z_{r,n}\!\ge\!1\,$ by construction. Note that 
 \bb
 \p(S_r\!>\!1) \!&=&\! 
 \p\Big( \cup_{1\le i<j\le r} \{X_{n,i}\!\ne\!0, X_{n,j}\!\ne\!0\}\Big)\\ 
 \!&\le&\! r\sum_{i=1}^{r-1} \p\(X_{n,i+1}\!\ne\!0, X_{n,1}\!\ne\!0\).  
 \ee 
   Thus, $\,\p(S_r\!>\!1) = o(r/n)\,$ if $\,(D')\,$ holds. In view of (\ref{13.6}), 
$\,\p(\z_{r,n}\!>\!1) \to 0\,$ as $\,n\to\infty$, i.e., (\ref{1}) holds.\\ 

  {\bf Remark 2.3}. If conditions $\,\Delta\,$ and $\,(D')\,$ hold, then (\ref{13.1}) is equivalent to $$ \lim_{n\to\infty}n\p(X_{n,1}\!\ne\!0) = \l. \eqno(\ref{GM-2}') $$ 
  Indeed, this follows from (\ref{13.6}), (\ref{rp}) and (\ref{2b2}) (cf. \cite{LLR}, Theorem 3.4.1).\\

	A generalisation of Corollary \ref{CorGM} to the case of stationary $\vp$-mixing r.v.s has been given by Utev \cite{Ut92}, Theorem 10.1, who has shown that conditions (\ref{GM-2}$'$) and ($D'$) are necessary and sufficient for (\ref{dep}). 
   Sufficient conditions for Poisson convergence without assuming stationarity have  been provided by Sevastyanov \cite{Sev}.   A Poisson limit theorem in the case of a two-dimentional random field $\,\{X_{i,j}\}\,$ has been given by Banis \cite{B85}.\\

   \pr of Theorem \ref{Dep}. Let $\,\{r\!=\!r_n\}\,$ be an arbitrary sequence from $\,{\cal R}$. Condition $\,\Delta\,$ and Lemma 2.4.1 from \cite{LLR} imply that for any $\,t\!\in\!\R,$ as $\,n\to\infty$, 
 \b                                                     \label{13.8} 
 \E \exp\(it S_n\) \!&=&\! \exp\!\( \frac{n}{r}\p(S_r\!>\!0) \E \left\{
 e^{itS_r}\!-1 | S_r\!>\!0 \right\} \) \,+\, o(1),\\ 
 \p(S_n=0) \!&=&\! \p^{n/r}(S_r\!=\!0) + o(1) \ =\  
 \exp\!\( -\frac{n}{r}\p(S_r\!>\!0) \) + o(1) 					\label{2b2} 
 \e 
(cf. (5.10) in \cite{N11}). 

  If (\ref{dep}) holds, then so does (\ref{13.1}): $\,\p(S_n\!=\!0)\to e^{-\l}\,$ as $\,n\!\to\!\infty$. Note that (\ref{13.1}) and (\ref{2b2}) yield (\ref{13.6}). 
		Since 
  $$ \E e^{itS_n} \to \exp(\l(e^{it}\!-\!1))\qquad (\forall t\!\in\!\R)\eqno(\ref{dep}^*)  $$ 
by the assumption, (\ref{13.8}) and (\ref{13.6}) entail $\,\E e^{it\z_{r,n}} \to e^{it},$ i.e., (\ref{1}) holds. 

   On the other hand, if (\ref{1}) and (\ref{13.1}) hold for some $\,\{r\!=\!r_n\}\!\in\!{\cal R}$, then (\ref{13.6}) is valid. Relations (\ref{13.6}) and (\ref{13.8}) yield (\ref{dep}$^*$). \hspace*{\fill} $\Box$\\

			\section{Accuracy of Poisson approximation}\label{2}

 	The problem of evaluating the accuracy of Poisson approximation to the distribution of a sum $$\,S_n=X_1\!+...+\!X_n\,$$ of independent 0-$\!$1 random variables has attracted a lot of attention among researchers (cf. \cite{BHJ,N11} and references wherein).   
		
  A natural task is to obtain a sharp estimate of the accuracy of Poisson approximation to the distribution of $\,\L(S_n)$. In this section we overview available estimates. 

	Historically, the accuracy of Poisson approximation was first studied in terms of the uniform distance (sometimes called the {\it Kolmogorov} distance). 

   The {\it uniform distance} $\,d_K(X;Y) \equiv d_K(F_X;F_Y)\,$ between the distributions of random variables $\,X\,$ and $\,Y\,$ with distribution functions (d.f.s) $\,F_X\,$ and $\,F_Y\,$ is defined as 
	$$ d_K(F_X;F_Y) = \sup_x|F_X(x)-F_Y(x)|. $$  

   Many authors evaluated the accuracy of Poisson approximation to $\,\L(S_n)\,$ in terms of the {\it total variation distance}. Recall that the total variation distance $\,\d(X;Y)\,$ between the distributions of r.v.s $\,X\,$ and $\,Y\,$ is defined as 
 \bb  
 \d(X;Y) \equiv \d(\L(X);\L(Y)) = 
 \sup_{A\in{\cal A}} \left|\p(X\!\in\!A)-\p(Y\!\in\!A)\right|, 
 \ee 
where $\,{\cal A}\,$ is a Borel $\sigma$-field. Evidently, $\,d_K(X;Y) \le \d(X;Y).$ 
   Note that $$ \d(X;Y) = \inf_{X',Y'}\p(X'\ne Y'), $$ 
where the infimum is taken over all random pairs $\,(X',Y')\,$ such that 
$\,\L(X')=\L(X)\,$ and $\,\L(Y')=\L(Y)\,$ \cite{Dob,B03}. 	
	
   The Gini-Kantorovich distance between the distributions of r.v.s $\,X\,$ and $\,Y\,$ with finite first moments (known also as the Kantorovich--Wasserstein distance) is 
 \beq                                                       \label{dV}
d_{_{G}}(X;Y)\equiv d_{_{G}}(\L(X);\L(Y))=\sup_{g\in{\cal F}}\left|\E g(X)-\E g(Y)\right|, 
 \eeq
where $\,{\cal F} = \{g\!: |g(x)\!-\!g(y)|\le|x\!-\!y|\}\,$ is the set of Lipschitz
functions. Note that 
 \beq                                                   \label{V} 
 d_{_{G}}\!\(X;Y\) = \inf_{X',Y'} \E|X'-Y'|,
 \eeq 
where the infimum is taken over all random pairs $\,(X',Y')\,$ such that 
$\,\L(X')=\L(X)\,$ and $\,\L(Y')=\L(Y)\,$ \cite{Val73}. 
	If $\,X\,$ and $\,Y\,$ take values in $\,\Z,$ then \cite{Pro56} 
  $$	d_{_{G}}\!\(X;Y\) = \sum_{i\ge 1}|\p(X\!\ge\!i)-\p(Y\!\ge\!i)|.	$$ 
	
   Distance $\,d_{_{G}}\,$ was introduced by Kantorovich \cite{Kan42} (to be precise, Kantorovich has introduced a class of distances that includes $\,d_{_{G}}$). We add the name of Gini since Gini \cite{Gi} used $\,\E|X\!-\!Y|$-type quantities. 
	Barbour et al. \cite{BHJ} called $\,d_{_{G}}\,$ the ``Wasserstein distance'' after Dobrushin \cite{Dob} attributed it to Vasershtein \cite{Va69}. 

  If distributions $\,P_1\,$ and $\,P_2\,$ have densities $\,f_1\,$ and $\,f_2\,$ with respect to a measure $\,\mu,$ set 
  $$ 
 d_{_H}^2(P_1;P_2) := \frac12 \int\!\(f_1^{1/2}-f_2^{1/2}\)^2 d\mu =
 1-\int\!\sqrt{f_1f_2} \,d\mu. 
 $$ 
Then $\,d_{_H}\,$ denotes the Hellinger distance. It is known that 
 \beq \label{d2} d_{_H}^2 \le \d \le d_{_H}\sqrt{2\!-\!d_{_H}^2}\,. \eeq 

  Denote 
	$$ \chi^2(P_1;P_2) = \int_{{\rm supp}P_2}\!\(dP_1/dP_2-1\)^2 dP_2. $$ 
  By the Cauchy-Bunyakovski  inequality, $$ 2\d(P_1;P_2) \le \chi(P_1;P_2). $$ 
	
  We denote by 
$$ d_{_{KL}}^2\!(P_1;P_2) \!= \int_{{\rm supp}P_2} \ln\!\(dP_1/dP_2\) dP_1 $$ 
the Kullback--Leibler divergence. 
  According to Pinsker's inequality,\index{Pinsker's inequality} 
  \beq \label{Pinsker} \d \le d_{_{KL}}/\sqrt{2}\,. \eeq 

  Though $\,d_{KL}^2\,$ is not a metric, it plays a role in statistics (cf. \cite{HR04}) and in the theory of large deviations (cf. \cite{N11}, p. 324, ex. 41). 
 
	Certain other distances can be found in \cite{Log90,N11,Roos01}. Below we present estimates of the accuracy of Poisson approximation for $\,\L(S_n)\,$ in terms of $\,d_K,$ $\,\d\,$ and $\,d_{_{G}}\,$ distances. 

					\subsection{Independent Bernoulli r.v.s}\label{2.1} 

   We denote by $\,{\bf B}(n,p)\,$ the Binomial distribution with parameters $\,n\,$ and $\,p$. Let $\,{\bf\Pi}(\l)\,$ denote the Poisson distribution with parameter $\,\l$; we denote by $\,\pi_\l\,$ a Poisson $\,{\bf\Pi}(\l)\,$ random variable. 
	
  Let $\,X_1,X_2,...,X_n\,$ be independent Bernoulli $\,{\bf B}(p_i)\,$ r.v.s. Denote $\,\l = \E S_n,$ 
	$$ p_i=\p(X_i\!=\!1)\ \ (i\!\ge\!1),\ \ 
	\l_k=\sum_{i=1}^n p_i^k\ \ \ (k\!\ge\!2),\ \  
	\theta = \sum^n_{i=1} p_i^2/\l. $$ 

  Many authors worked on the problem of evaluating the accuracy of Poisson approximation to $\,\L(S_n)\,$ in terms of the uniform distance $\,d_K\,,$ the total variation distance $\,\d\,$ and the Gini--Kantorovich distance $\,d_{_{G}}$. 

   It seems natural to approximate $\,{\bf B}(n,p)\,$ by the Poisson distribution. For instance, in the case of identically distributed Bernoulli $\,{\bf B}(p)\,$ r.v.s $\,\{X_i\}\,$ one has 
	\beq \label{Usl}
	\p(S_n\!=\!k) \equiv \p(S_n\!=\!k| N_n\!=\!n) = 
	\p(\pi_n(p)\!=\!k|\pi_n(1)\!=\!n), 
	\eeq
where $\,N_n\!\equiv\!n\,$	is the total number of 0's and 1's among $\,X_1,X_2,...,X_n\,$ and $\,\{\pi_n(t),t\!\in\![0;1]\}\,$ is a Poisson jump process on $\,[0;1]\,$ with intensity rate $\,n$. Thus, 
  $$ {\bf B}(n,p) = \L(\pi_n(p)|\pi_n(1)\!=\!n). \eqno(\ref{Usl}^*) $$ 
	
	   Tsaregradskii \cite{Ts58} has shown that 
 \beq 																	\label{Ts}
 d_K(F_X;F_Y) \le \int_{-\pi}^\pi \frac{|\E e^{itX}-\E e^{itY}|}{4|t|}dt 
 \eeq 
if $\,X\,$ and $\,Y\,$ are integer-valued r.v.s, and derived the estimate 
 \beq \label{Tsa58} 
 d_K({\bf B}(n,p);{\bf\Pi}(np)) \le p\pi^2e^{2p(2-p)}/16(1\!-\!p) 
 \qquad(p\!\in\!(0;1/2]). \eeq 
Note that $\,\pi^2/16 \approx 0.617$. 
   Inequality (\ref{Tsa58}) seems to be the first estimate of the accuracy of Poisson approximation with explicit constant. 
	
	In the case of non-identically distributed Bernoulli $\,{\bf B}(p_i)\,$ random variables Franken \cite{Fr} has shown that 
	$$ d_K(S_n;\pi_\l) \le 0.6p^*_n $$ if $\,p^*_n := \max_{i\le n}p_i \le 1/4$. 
  Shorgin\index{Shorgin} \cite{Sho77} has proved that 
 $$ d_K(S_n;\pi_\l) \le c_1\theta/(1\!-\!\sqrt{\theta}\,) \qquad(\theta\!<\!1) $$ 
where $\,c_1 = (1\!+\!\sqrt{\pi/2})/2 < 1.13$.  
  According to Daley \& Vere-Jones \cite{DV08}, $$ d_K(S_n;\pi_\l) \le 0.36\theta. $$ 
  Roos \cite{Roos01} has shown that 
	$$ d_K(S_n;\pi_\l) \le \Big(1/2e + 1.2\sqrt{\theta}/(1\!-\!\sqrt{\theta})\Big)\theta. $$ Note that $\,1/2e \le 0.184$.\\ 

   Kontoyiannis et al. \cite{KHJ}\index{Kontoyiannis} have shown that 
 $$ d_{\!_H}^2\!(S_n;\pi_\l) \le \l^{-1}\sum^n_{i=1} p_i^3/(1\!-\!p_i). $$
   Borisov \& Vorozheikin \cite{BV} present sharp lower and upper bounds to $\,\chi^2\!({\bf B}(n,p);{\bf\Pi}(np))\!:$ 
  $$ 0\le \chi^2\!({\bf B}(n,p);{\bf\Pi}(np)) - p^2/2 - 2p^3/3n \le 
	p^4/(1\!-\!p) + p^8(23\!-\!20p)/(1\!-\!p)^2. $$ 
  Harremo\"es \& Ruzankin \cite{HR04} present lower and upper bounds to $\,d_{_{KL}}^2\!({\bf B}(n,p);{\bf\Pi}(np)).$ In particular, they have shown that 
  $$ 2d_{_{KL}}^2\!({\bf B}(n,p);{\bf\Pi}(np)) = 
	(-p-\ln(1\!-\!p))(1\!+\!O(1/n)). $$ $\,$

	 Many authors worked on the problem of evaluating the total variation distance $\,\d\!\(S_n;\pi_\l\)\,$ (cf. \cite{BHJ,N11} and references wherein). 
		Prohorov \cite{Pro53} has established the existence of an absolute constant $\,c\,$ such that 
 \beq                                   \label{Proh} 
 \d({\bf B}(n,p);{\bf\Pi}(np)) \le cp. \eeq 
	Kolmogorov \cite{K56} points out that 
	$$ \d(S_n;\pi_\l) \le C\sum^n_{i=1} p_i^2, $$ 
where $\,C\,$ is an absolute constant. LeCam \cite{LeCam60,LeCam65} attributes inequality 
	\beq   \label{LeCam} 
	\d(S_n;\pi_\l) \le \sum^n_{i=1} p_i^2 
	\eeq 
to Khintchin \cite{Khi}. Bound (\ref{LeCam}) is sharp: according to (2.10) in Deheuvels \& Pfeifer \cite{DP88}, $$\,\d(S_n;\pi_\l) \ge np^2(1\!+\!O(p))\,$$ 
in the case of i.i.d. Bernoulli $\,{\bf B}(p)\,$ r.v.s if	$\,np\to0$. 

   Note that (\ref{LeCam}) is a consequence of the property of $\,\d\,$ and the following fact: 
	\beq \label{TVD} \d({\bf B}(p);{\bf\Pi}(p)) = (1\!-\!e^{-p})p \le p^2. 
	\eeq 
  Indeed, denote $\,\bar X\!=\!(X_1,...,X_n),$ $\,\bar\pi\!=\!(\pi_{p_1},...,\pi_{p_n}),$ where $\,\{\pi_{p_i}\}\,$ are independent Poisson $\,{\bf\Pi}(p_i)\,$ r.v.s. Then 
  $$ \d(S_n;\pi_\l) \le \d(\bar X;\bar\pi) \le \sum^n_{i=1} 
	   \d(X_i;\pi_{p_i}) \le \sum^n_{i=1} p_i^2. \eqno(\ref{LeCam}') $$ 

   Set $\,\tilde p_i = -\ln(1\!-\!p_i)\ \ (1\!\le\!i\!\le\!n),$ and put $\,\mu=\sum^n_{i=1} \tilde p_i$. According to Serfling \cite{S75}, 
	$$ \d(S_n;\pi_{\mu}) \le \sum^n_{i=1} \tilde p_i^2/2. $$ 
	Kerstan \cite{K64} has shown that 
	\beq \label{K64} \d(S_n;\pi_\l) \le 1.05\theta. \eeq 
Romanowska \cite{Ro77} has noticed that 
 \beq  \label{Rom} 
 \d({\bf B}(n,p);{\bf\Pi}(np)) \le p/2\sqrt{1\!-\!p}\,. \eeq 
	The popular estimate 
 \beq                                              \label{be}
 \d(S_n;\pi_\l) \le \l^{-1}(1\!-\!e^{-\l}) \sum^n_{i=1} p_i^2 
 \eeq 
is effectively due to Barbour and Eagleson \cite{BE}. 

   Presman \cite{Pre85} has established an estimate of $\,\d(S_n;\pi_\l)\,$ 
with the constant $0.83$ at the leading term. 
   In the case of i.i.d. Bernoulli $\,{\bf B}(p)\,$ r.v.s Presman's bound becomes 
 \beq \label{Pre85} 
 \d({\bf B}(n,p);{\bf\Pi}(np)) \le 0.83 p/(1\!-\!p)(1\!-\!1/n). \eeq 
Xia \cite{Xia97} has derived an estimate with the constant $0.6844$ at the leading term.  

	Roos \cite{Roos01} (see also $\rm\check{C}$ekanavi\v{c}ius \& Roos \cite{CR06}) has obtained a bound with a correct constant $\,3/4e\approx0.276\,$ at the leading term: if $\,\theta\!<\!1,$ then 
 \beq                                                  \label{R2001}
 \d\!\(S_n;\pi_\l\) \le 3\theta/4e (1\!-\!\sqrt{\theta}\,)^{3/2}\,. 
 \eeq 
  Note that $\,\theta/(1\!-\!\sqrt{\theta}\,)^{3/2} \ge \theta(1\!+\!1.5\sqrt{\theta}+\!3.75\theta).$ 
	
  Roos \cite{Roos01} has shown also that $$ \d(S_n;\pi_\l) \sim 3\theta/4e $$ 
if $\,\theta\!\to\!0\,$ and $\,\l\!\to\!1\,$ as $\,n\!\to\!\infty$.	
Thus, constant $\,3/4e\,$ cannot be improved. 

   Denote 
	\bb 
	p^*_n &=& \max_{i\le n}p_i\,,\qquad\qquad\!\ \ve = \min\!\left\{1; 
	\(2\pi[\l\!-\!p^*_n]\)^{-1/2} + 2\delta/(1\!-\!p^*_n/\l)\right\},\\ 
	\delta &=& \frac{1\!-\!e^{-\l}}\l\sum^n_{i=1} p_i^2\,,\,\ \ \ 
  \delta^* = \frac{1\!-\!e^{-\l}} \l\sum^n_{i=1} p_i^3 \,.	
	\ee	
Note that $\,\delta^2\le\delta^*\,.$ 
   The following inequality from \cite{N11}, Theorem 4.12, sharpens the second-order term of the right-hand side (r.h.s.) of estimate (\ref{R2001}): 
 \beq                                     	\label{RN}  
 \d(S_n;\pi_\l) \le 3\theta/4e + 2\delta^*\ve + 2\delta^2. 
 \eeq 
   In the case of $\,\L(S_n)={\bf B}(n,p)\,$ estimate (\ref{RN}) becomes 
  $$ 
   \d(S_n;\pi_{np}) \le 3p\!/4e + 2(1\!-\!e^{-np})p^2\ve 
   + 2(1\!-\!e^{-np})^2p^2, \eqno(\ref{RN}^*) 
	$$ 
where $\,\ve = \min\{1; \(2\pi[(n\!-\!1)p]\)^{-1/2}\!+ 2(1\!-\!e^{-np})p/(1\!-\!1/n)\}.$ The second-order term in (\ref{RN}$^*$) is of order $\,p^2 \wedge np^3$.\\ 
	
	In applications one often has $\,\l\equiv\l(n)\to\infty\,$ as $\,n\to\infty.$ 
Hence estimates with the ``magic factor'' $\,(1\!-\!e^{-\l})/\l\,$ attract special interest. 

  The possibility of the ``super--magic'' factor $\,e^{-\l}\,$ when one approximates $\,S_n\!\in\!A\,$ for a bounded $\,A\,$ has been discussed in \cite{N11}, ch. 4.5 (such approximations are of interest in extreme value theory). 
	For instance, if $\,\{X_i\}\,$ are independent Bernoulli $\,{\bf B}(p_i)\,$ r.v.s and $\,A=\{0\}$, then 
	\beq \label{supe} 0\le \p(\pi_\l\!=\!0) - \p(S_n\!=\!0) \le 
	e^{-\l+p^*_n} \sum^n_{i=1} p_i^2\!/2. \eeq 
Indeed, set $\,c_i=\prod^{i-1}_{j=1}e^{-p_j}\!\prod^n_{j=i+1}(1\!-\!p_j).$ 
Since $\,e^{-p_i}\!-\!1\!+\!p_i \le p_i^2/2\,$ by Taylor's formula, 
  $$ \p(\pi_\l\!=\!0)-\p(S_n\!=\!0) = \sum^n_{i=1} (e^{-p_i}\!-\!1\!+\!p_i)c_i 
	   \le e^{-\l+p^*_n} \sum^n_{i=1} p_i^2\!/2. $$ 
  In the case of the Binomial $\,{\bf B}(n,p)\,$ distribution (\ref{supe}) becomes 
	$$ 0\le \p(\pi_{np}\!=\!0) - \p(S_n\!=\!0) \le \frac12 np^2e^{-(n-1)p} 
					  \le \frac{2n}{e^2(n\!-\!1)^2}\,. \eqno(\ref{supe}^*) $$ 
Note that $\,2/e^2\approx0.2707$. 
	
  Bound (\ref{RN}) is a consequence of inequality 
 $$ |\p(S_n\!\in\!A)-\p(\pi_\l\!\in\!A)| \le |\p(\pi_\l\!+\!1\!\in\!A) - 
 \p(\pi_\l^\star\!\in\!A)|\theta/2 +2\delta^*\ve+2\delta^2. \eqno(\ref{RN}^\star) $$ 
The first term on the r.h.s. of (\ref{RN}$^\star$) has the ``super--magic'' factor if $\,A\,$ is finite.\\ 
	
   Estimates in terms of the Gini-Kantorovich distance are available as well. 
Denote $\,\mu = -\sum_{i=1}^n\ln(1\!-\!p_i).$ If $\,p^*_n\le1/2$, then 
	\beq \label{DF89} d_{_{G}}(S_n;\pi_\mu) \le \sum_{i=1}^n p_i^2/2(1\!-\!p_i) 
  \eeq (Deheuvels et al. \cite{DF89}).  
  Witte \cite{W90} has shown that 
	\beq \label{W14} 
	d_{_{G}}(S_n;\pi_\l) \le \frac{-\sqrt{e\l}}{2\sqrt{2\pi}} 
	\ln(1\!-\!2\theta e^{2p_n^*}). 
	\eeq  
According to \cite{N11}, formula (4.53), 
  \beq \label{2b-1}  d_{_{G}}(S_n;\pi_\l) \le 
  \Big(1\wedge\frac{_4}{^3}\sqrt{2/e\l}\,\Big) \sum_{i=1}^n p_i^2 . 
  \eeq  
Roos \cite{Roos01} has shown that 
  \beq \label{R01} d_{_{G}}(S_n;\pi_\l) \le \(1/\sqrt{2e} + 
	   1.6\sqrt{\theta}(2\!-\!\theta)/(1\!-\!\sqrt{\theta})\) \theta\sqrt{\l}\,. \eeq A recent survey is Zacharovas \& Hwang \cite{ZH08}.\\ 

  Sharp non-Poisson approximation to the Binomial $\,{\bf B}(n,p)\,$ distribution function $\,\p(S_n\!\le\!\cdot)\,$ has been given by Zubkov \& Serov \cite{ZS}. 
	
	Denote by $\,\Phi\,$ the standard normal d.f.. Let 
	$$ \Lambda(x) = x\ln(x/p) + (1\!-\!x)\ln\((1\!-\!x)/(1\!-\!p)\) \qquad(0\!<\!x\!<\!1) $$ denote the rate function of the Bernoulli distribution $\,{\bf B}(p)\,$ (cf. \cite{N11}, p. 322), and set 
  $$ Z_{n,p}(k) = \Phi\Big(\hbox{sgn}(k/n\!-\!p) \sqrt{2n\Lambda(k/n)}\,\Big). $$ 
  Then \cite{ZS} 
	\beq \label{ZS} Z_{n,p}(k) \le \p(S_n\!\le\!k) \le Z_{n,p}(k\!+\!1). \eeq 
	
  The following large deviations inequality is due to Bernstein \cite{B46}, p. 168: $$ \p((S_n\!-\!np)/\sqrt{npq} > t\(1+\gamma_{p,n}\) < \exp(-t^2/2) \qquad(t\!>\!0), $$ where $\,q=1\!-\!p,$ $\,\gamma_{p,n} = t_{p,n}(q\!-\!p)/6 +	t_{p,n}^2(p^3\!+\!q^3)/12,$ $\,t_{p,n} = t/\sqrt{npq}\,$.\\

			{\bf Asymptotics of} {\boldmath$\,\d(S_n;\pi_\l)$}. 
	The asymptotics of $\,\d(S_n;\pi_\l)\,$ in the case of identically distributed Bernoulli $\,{\bf B}(p)\,$ r.v.s has been established by Prohorov \cite{Pro53}: 
 \beq                                         \label{Pro} 
 \d({\bf B}(n,p);{\bf\Pi}(np)) = p/\sqrt{2\pi e} 
 \Big(1+O\big(1\!\wedge\!(p+\!1/\!\sqrt{np}\,)\big)\Big). 
 \eeq 
Kerstan \cite{K64}, Deheuvels \& Pfeifer \cite{DP86,DP88AISM}, Deheuvels et al. \cite{DF89} and Roos \cite{Roos99} have generalised (\ref{Pro}) to the case of non-identically distributed 0-$\!$1 r.v.s.  
   Deheuvels \& Pfeifer \cite{DP86} present also the asymptotics of $\,\d(S_n;\pi_\l)\,$ in the case where $\,\l\!\to\!{\rm const}\,$ as $\,n\!\to\!\infty$. 
	
	The following result concerning the asymptotics of $\,\d(S_n;\pi_\l)\,$ uses the notation from \cite{N11}, ch. 4. 
		Given a non-negative integer-valued random variable $\,Y,$ we denote by $\,Y^\star\,$ a random variable with the distribution 
 \beq \label{star} 
 \p(Y^\star\!=\!k) = \p(Y\!=\!k)(k\!-\!\l)^2/\l \qquad (k\!\in\!\Z\!). \eeq 
The next bound is a consequence of Theorem \ref{asy}.

  \begin{theorem} If $\,X_1,...,X_n\,$ are independent Bernoulli r.v.s, 
$\,\L(X_i)={\bf B}(p_i),$ then 
 \beq                                                           \label{asTVD}
 |\d(S_n;\pi_\l) - \theta\d(\pi_\l^\star;\pi_\l\!+\!1)/2| \le 2\delta^*\ve+2\delta^2. \eeq \end{theorem}

   One can check that 
 \beq                                       \label{TVDs}
 \d(\pi_\l^\star;\pi_\l\!+\!1) = \sqrt{2/\pi e}+O\big(1/\sqrt{\l}\,\big) \eeq	
as $\,\l\to\infty$. Thus, 
 \beq                                                        \label{aTVD}
 \d(S_n;\pi_\l) = \theta/\sqrt{2\pi e} \(1+O\big(\theta\!+\!1/\sqrt{\l}\,\big)\) 
 \eeq if $\,\l\to\infty\,$ and $\,\theta\to0\,$ as $\,n\!\to\!\infty$.\\

   \noindent{\bf Example 2.1}. Let $\,p_i\!=\!1/i,$ $\,i\!\in\!\N.$ Then $\,p^*_n\!=\!1,$ $\,\l=\l(n)\to\infty,$ $\,\theta\to0\,$ as $\,n\!\to\!\infty$, and (\ref{aTVD}) entails $\,\d(S_n;\pi_\l) \sim \theta/\sqrt{2\pi e}\,.$ \hspace*{\fill}$\Box$\\ 

   Deheuvels et al. \cite{DF89} have shown that 
  $$ \big|d_{_{G}}(S_n;\pi_\l) - \l_2e^{-\l}\l^{[\l]}/[\l]!\big| \le 
	   2(2\theta)^{3/2}\sqrt{\l}/(1\!-\!\sqrt{2\theta}) \qquad(\theta\!<\!1/2). $$ 
  Borisov \& Vorozheikin \cite{BV} present asymptotic expansions of 
$\chi^2\!({\bf B}(n,p);{\bf\Pi}(np))$. 
		
		$\,$

   {\bf Shifted Poisson approximation}. 
   Shifted (translated) Poisson approximation to $\,{\bf B}(n,p)\,$ has been considered by a number of authors (see \cite{BC02,BX06,CV01,K86,N18} and references therein). 
	 The accuracy of shifted Poisson approximation can be sharper than that of pure Poisson approximation. 
	 Another advantage of using shifted Poisson approximation is the possibility to derive a more general result (e.g., a	uniform in $\,p\,$ estimate of $\,\d({\bf B}(n,p);{\bf\Pi}(np))$, cf. (\ref{PEx3}) below). 
  
  Let $\,X_1,...,X_n\,$ be independent 0-$\!$1 r.v.s. 
	Set $\,p_i\!=\!\p(X_i\!=\!1),$ $\,q_i\!=\!1\!-\!p_i,$ 
	 $$ \l=\E S_n,\ \sigma^2\!=\!\hbox{var}\,S_n,\ \l_2=\l\!-\!\sigma^2. $$ 
  Denote $\,[x] = \max\{k\!\in\!{\bf Z}\!: k\!\le\!x\},$ $\,\{x\} = x\!-\![x].$ 	
	We define r.v. $$ Y = [\l_2]\!+\!\pi_{\l-[\l_2]}. $$ 
	Note that $\,\hbox{var}\,\pi_{\l}\!-\!\hbox{var}\,S_n \!=\! \l_2,$ while 
$\,\hbox{var}\,Y\!-\!\hbox{var}\,S_n =\{\l\} \!<\!1.$ 

  The following result is due to $\rm\check{C}$ekanavi\v{c}ius \& Vaitkus \cite{CV01}. 

	\begin{theorem} \label{CV2001} 
If $\,\sigma^2\!\ge\!4,$ then 
	\beq  \label{CV} 
	\d(S_n;Y) \le 0.93 \sigma^{-3}\l_2  
	+ \{\l_2\}/(\sigma^2\!+\!\{\l_2\}) + e^{-\sigma^2/4} . 
	\eeq \end{theorem}

  Let $\,\{X_i\}\,$ be i.i.d. Bernoulli $\,{\bf B}(p)\,$ r.v.s. Then the right-hand side (r.-h.s.) of (\ref{CV}) is $$O\Big(\sqrt{p/n}+\!1/np\Big).$$ A similar bound in terms of the uniform distance has been established by Kruopis \cite{K86}. 
	
  Set $\,q=1\!-\!p$, where $\,0\!<\!p\!<\!1.$ Then 
	$$ Y = [np^2]\!+\!\pi_{npq+\{np^2\}},\ \ \E Y = np,\ \ 
	   \hbox{var} Y = npq\!+\!\{np^2\}.	$$ 
The following Theorem \ref{T2018} presents a uniform in $\,p\!\in\![0;1/2]\,$ bound to $\,\d(S_n;Y)$.

	\begin{theorem} \label{T2018} {\rm\cite{N18}} As $\,n\!>\!4$, 
 \beq  \label{PEx3} 
 \sup_{0\le p\le1/2} \d(S_n;Y) \le 
 \frac2{\sqrt{\pi}} \frac{1\!+\!1/\sqrt{e}}{\sqrt{n}-\!2} 
	+ \frac{2\!+\!4/\sqrt{\pi}}{n\!-\!2\sqrt{n}}\ \qquad(n\!>\!4).
 \eeq 	\end{theorem}

  Theorem \ref{T2018} can be compared with the Berry--Esseen inequality 
	$$ d_K({\bf B}(n,p);{\cal N}(np,npq)) \le C/\sqrt{np}\, $$ 
(see, e.g., \cite{Shev11}) as well as with the results by Meshalkin \cite{Me60} and Pressman \cite{Pre83}. 
  Estimate (\ref{PEx3}) is uniform in $\,p\!\in\![0;1/2]$. Note that a uniform in $\,p\!\in\![0;\!1/2]\,$ Berry--Esseen estimate would be infinite. 
   Inequality (\ref{PEx3}) has advantages over Meshalkin's \cite{Me60} and Pressman's \cite{Pre83} results as the constants in (\ref{PEx3}) are explicit (which matters in applications); besides, the structure of the approximating distribution $\,\L(Y)\,$ is simpler and does not assign mass to negative numbers.   Bound (\ref{PEx3}) is preferable to (\ref{be}) -- (\ref{RN}) if $\,p\!>\!4e/\sqrt{n}\,.$\\ 

	An estimate of the accuracy of shifted Poisson approximation to the distribution of a sum of Bernoulli $\,{\bf B}(p_i)\,$ r.v.s in terms of the Gini-Kantorovich distance has been given by Barbour \& Xia \cite{BX06} in the assumption that $\,\l_2\,$ is an integer.\\

				{\bf Poisson approximation to the multinomial distribution}. 
Results on the accuracy of Poisson approximation to the distribution of a sum of Bernoulli r.v.s can be generalised to the case of a multinomial distribution. 

   Let $\,\bar S_n\,$ be a random vector with multinomial distribution $\,{\bf B}(n,p_1,...,p_m)$: 
 \beq                                                   \label{2-43}
 \p(\bar S_n = \bar l\,) = \frac{n!} {l_1!...l_m!(n-l)!}
 p_1^{l_1} ... p_m^{l_m} (1\!-\!p)^{n-l} \,,
 \eeq 
where $\,l_i\!\in\!\Z\ (\forall i),$ $\,\bar l = (l_1,..,l_m),$ 
$\,l=l_1+\!...\!+l_m\le n,$ $\,p=p_1+\!...\!+p_m.$ 

   Formula (\ref{2-43}) describes, in particular, the joint distribution of the increments of the empirical d.f.. 

   Note that 
 \beq                                   \label{2-sum}
 \bar S_n \st \bar\xi_1+\!...\!+\bar\xi_n,
 \eeq
where $\,\bar\xi, \bar\xi_1,..., \bar\xi_n\,$ are i.i.d. random vectors with the distribution 
 $$
 \p(\bar\xi\!=\!\bar0)=1\!-\!p\,,\ 
 \p(\bar\xi\!=\!\bar e_j) = p_j \qquad(1\!\le\!j\!\le\!m), 
 $$ 
vector $\,\bar e_j\,$ has the $\,j^{\rm th}\,$ coordinate equal to 1 and the other coordinates equal to 0.

   Let $$ \bar\pi = \(\pi_1,...,\pi_m\) $$ 
be a vector of independent Poisson r.v.s with parameters $\,np_1,...,np_m,$ and let $\,\pi_n(\cdot)\,$ denote a Poisson jump process on $\,[0;1]\,$ with intensity rate $\,n$. Then $\,\bar\pi\,$ is a vector of increments of process $\,\pi_n(\cdot)$: 
$\,\pi_1\st\pi_n(p_1)$,...,$\,\pi_m\st\pi_n(p)-\pi_n(p\!-\!p_m)$. 
  Note that 
	$$ \p(\bar S_n = \bar l\,) = \p\(\pi_n(p_1) \!=\! l_1,..., 
	\pi_n(p)-\pi_n(p\!-\!p_m) \!=\! l_n \,|\,\pi_n(1)\!=\!1\) \eqno(\ref{2-43}^*)$$ (cf. (\ref{Usl})). 

    Arenbaev \cite{Ar76} has shown that 
 \beq \label{Aren} 
 \d(\bar S_n;\bar\pi) = p/\sqrt{2\pi e}\(1+O(1\!\wedge\!1/\!\sqrt{np}\,)\) 
 \eeq 
if $\,n\!\to\!\infty\,$ (the term $\,1/\!\sqrt{np}\,$ in (\ref{Aren}) apparently needs to be replaced with $\,p+\!1/\!\sqrt{np}\,,$ cf. (\ref{Pro})). 
   Arenbaev (\cite{Ar76}, formulas (5)--(9$^\prime$)) has shown also that 
 \beq                                                   \label{tdvM}
 \d(\bar S_n;\bar\pi) = \d({\bf B}(n,p);{\bf \Pi}(np)). 
 \eeq 
Using (\ref{tdvM}) and (\ref{RN}), we deduce 
 \beq \label{72} 
 \d(\bar S_n;\bar\pi) \le 3p/4e + 4(1-e^{-np})p^2. 
 \eeq 
   According to Deheuvels \& Pfeifer \cite{DP88}, 
 \beq \label{MultDP} 
 | \d(\bar S_n;\bar\pi) - K_{n,\l} | \le \max\{16p^2;5np^3\},
 \eeq 
where  $\,K_{n,\l} = np^2 e^{-np} \( (np)^{\a-np}(\a\!-\!np)/\a! - 
(np)^{\beta-np} (\beta\!-\!np)/\beta! \)\!/2,$ 
$$\,\a = np\!+\!1/2+\sqrt{np\!+\!1/4}\,,\ \,\beta = np\!+\!1/2-\sqrt{np\!+\!1/4}\,.$$ 

   The case of non-identically distributed random vectors $\,\bar\xi_1,..., \bar\xi_n\,$ has been treated by Roos \cite{R15}.  
   A generalisation of (\ref{72}) to the case of a stationary sequence of {\it dependent} r.v.s is given in \cite{N11}, Theorem 6.8.\\ 

   \noindent{\it Open problem}.\\ 
	 \noindent2.1. Improve the constants in (\ref{DF89})--(\ref{2b-1}).\\ 
	 \noindent2.2. Generalise Theorem \ref{T2018} to the case of $m$-dependent r.v.s. 

						\subsection{Dependent Bernoulli r.v.s}\label{2.2}

	We present below generalisations of (\ref{LeCam}) and (\ref{be}) to the case of dependent Bernoulli r.v.s. 

  Let $\,X_1,...,X_n\,$ be (possibly dependent) Bernoulli r.v.s. 
	Chen \cite{Chen75} pioneered the use of Stein's method in deriving estimate of the accuracy of Poisson approximation, and obtained an estimate of the accuracy of Poisson approximation to the distribution of a sum of $\,\vp$-mixing r.v.s. 
	
	Set $\,p_i\!=\!\p(X_i\!=\!1|X_1,...,X_{i-1}).$ 
	A generalisation of (\ref{LeCam}) has been given by Serfling \cite{S75}: 
 $$ 
 \d(S_n;\pi_\l) \le \sum^n_{i=1} (\E p_i)^2 + \sum^n_{i=1}\E|p_i\!-\!\E p_i|, 
 \eqno(\ref{LeCam}^*) $$ 
 \beq \label{S75} \ \  d_K(S_n;\pi_\l) \le 
 \frac2\pi\sum^n_{i=1}(\E p_i)^2 + \sum^n_{i=1}\E|p_i\!-\!\E p_i|. \eeq 
		
   Let $\,\{X_a, a\!\in\!J\}\,$ be a family of dependent Bernoulli $\,{\bf B}(p_a)\,$ random variables. Assign to each $\,a\!\in\!J\,$ a ``neighborhood'' $\,B_a\!\subset\!J\,$ such that $\,\{X_b, b\!\in\!J\!\setminus\!B_a\}\,$ are ``almost independent'' of $\,X_a\,$ (for instance, if $\,\{X_b\}\,$ are $m$--dependent r.v.s and $\,J=\{1,...,n\},$ then $\,B_a=[a\!-\!m;a\!+\!m]\cap J$). 
	
   The idea of splitting the sample into ``strongly dependent'' and ``almost independent'' parts goes back to Bernstein \cite{B26} (see also \cite{Sev}). 
	
	Denote $$ S = \sum_{a\in J} X_a,\ \l=\E S, $$ and let 
 \begin{eqnarray*}
 \delta_1 &=&\sum_{a\in J}\sum_{b\in B_a} \!\E X_a\E X_b\,,\ 
 \delta_2 =\sum_{a\in J}\sum_{b\in B_a\backslash \{a\}} \E X_a X_b,\\ 
 \delta_3 &=&\sum_{a\in J} \E\,\Big| \E X_a - \E 
 \Big\{ X_a \Big| \sum\nolimits_{b\in J\backslash B_a} \!X_b \Big\} \Big|. 
 \end{eqnarray*} 
   The following Theorem \ref{T37} is cited from Arratia et al. \cite{AGG} and Smith \cite{S88}.

    \begin{theorem} \label{T37} There holds 
 \beq                                                           \label{311}
 \d(S;\pi_\l) \le \frac{1\!-\!e^{-\l}}{\l} \Big(\delta_1 + \delta_2\Big) +
 \min\{1;\sqrt{2/e\l}\,\} \delta_3\,. \eeq \end{theorem}

   In the case of independent random variables one can choose $\,B_a=\{a\},$ 
then (\ref{311}) coincides with (\ref{be}). 

   Theorem \ref{T37} has applications to the problem of Poisson approximation to the distribution of the number of long head runs in a sequence of Bernoulli r.v.s, and to the problem of Poisson approximation to the distribution of the number of long match patterns in two sequences (e.g., DNA sequences, see \cite{BHJ,N11} and references therein). 

  The topic concerning $\,\L(S_n)\,$ in the case of stationary dependent r.v.s $\,\{X_i\}\,$ has applications in extreme value theory \cite{LLR,N11}. The case where the sequence $\,X_1,...,X_n\,$ is a moving average is related to the topic concerning the so-called Erd\"os--R\'enyi maximum of partial sums (cf. \cite{N11}, ch. 2). 

  Estimates of the accuracy of Poisson approximation for some special types of dependence among $\,\{X_a, a\!\in\!J\}\,$ can be found in Barbour et al. \cite{BHJ}. 	An estimate of the accuracy of shifted Poisson approximation to the distribution of a sum of dependent Bernoulli $\,{\bf B}(p_i)\,$ r.v.s in terms of the total variation distance is given by $\rm\check{C}$ekanavi\v{c}ius \& Vaitkus \cite{CV01}. 
  A generalization of Theorem \ref{T37} to the case of compound Poisson approximation has been given by Roos \cite{Roos94}.\\

   \noindent{\it Open problem}.\\ 
   \noindent2.3. Improve the constants in (\ref{311}).

					\subsection{Independent integer-valued r.v.s}\label{2.3}

   The topic of Poisson approximation to the distribution of a sum of integer-valued r.v.s has applications in extreme value theory, insurance, reliability theory, etc. (cf. \cite{BalK,BHJ,LLR,N11}). For instance, in insurance applications the sum $\,S_n=\sum_{i=1}^n\! Y_i\1\{Y_i\!>\!y_i\}\,$ of integer-valued r.v.s allows to account for the total loss from the claims exceeding excesses $\,\{y_i\}$. One would be interested if Poisson approximation to $\,\L(S_n)\,$ is applicable. 
	
   In extreme value theory one often deals with the number of extreme (rare) events represented by a sum $\,S_n=\xi_1+...+\xi_n\,$ of 0-$\!$1 r.v.s (indicators of rare events). The r.v.s $\,\xi_1,...,\xi_n\,$ can be dependent. 
   One way to cope with dependence is to split the sample into blocks, which can be considered almost independent (the so-called Bernstein's blocks approach \cite{B26}). 
	The number of r.v.s in a block is an integer-valued r.v.; thus, the number of rare events is a sum of almost independent integer-valued r.v.s. 
 
	In all such situations one deals with a sum of non-negative integer-valued r.v.s that are non-zero with small probabilities, and Poisson or compound Poisson approximation to $\,\L(S_n)\,$ appears plausible. 
	An estimate of the accuracy of Poisson approximation to the distribution of $\,S_n\,$ can indicate whether Poisson approximation is applicable. 
	
   The problem of evaluating the accuracy of Poisson approximation to the distribution of a sum of independent non-negative integer-valued r.v.s has been considered, e.g., in \cite{B87,BJ89,N11}.    Inequality (\ref{LeCam}) and the Barbour-Eagleson estimate (\ref{be}) have been  generalised to the case of non-negative integer-valued r.v.s by Barbour \cite{B87}. Theorem \ref{P38} below presents another result of that kind (see \cite{N11}, ch. 4.4). 
	
  Let $\,X_1,X_2,...,X_n\,$ be independent non-negative integer-valued r.v.s, 
	$$ S_n = X_1+...+X_n,\ \l = \E S_n, $$ 
$\,\pi_\l\,$ denotes a Poisson $\,{\bf\Pi}(\l)\,$ r.v.. 

	Franken \cite{Fr} has shown that 
	$$ d_K(S_n;\pi_\l) \le \frac2\pi \sum^n_{i=1} (\E^2 X_i+\E X_i(X_i\!-\!1)) 
	$$ 	
Denote $\,\l^* = \sum^n_{i=1}\!\p(X_i\!=\!1),$ $\,\l^*_2 = \sum^n_{i=1}\!\p(X_i\!=\!1)^2.$ 
  Kerstan \cite{K64} has proved that 
  $$ \d(S_n;\pi_{\l^*}) \le \sum^n_{i=1} \p(X_i\!\ge\!2) + 
	\min\{ \l^*_2 ; 1.05\l^*_2/\l^* \}. $$ 
An early survey on the topic is Witte \cite{W90}. 

  Given a random variable $\,Y\,$ that takes values in $\,\Z$, let $\,Y^*\,$ denote a random variable with the distribution 
 \beq                                                   \label{313}
 \p(Y^*\!=\!m)=(m\!+\!1)\p(Y=m\!+\!1)/\E Y \qquad(m\!\ge\!0). \eeq 
   Distribution (\ref{313}) differs by a shift from the distribution introduced by Stein \cite{St92}, p. 171. 
	 Note that $\,Y^*\st Y\,$ if and only if $\,\L(Y)\,$ is Poisson.

 	\begin{theorem} \label{P38} As $\,n\!\ge\!1,$ 
 \b                                               \label{314}
 \d(S_n;\pi_\l) &\le& \l^{-1}(1\!-\!e^{-\l}) \sum^n_{i=1} d_{_G}(X_i;X_i^*)\E X_i\,,\\
 d_{_G}(S_n;\pi_\l) &\le& \min\!\Big\{1;\frac{_4}{^3}\sqrt{2/e\l}\,\Big\}
 \sum^n_{i=1} d_{_G}(X_i;X_i^*)\E X_i\,. 					\label{415G} 
 \e \end{theorem}

	 In the case of Bernoulli $\,{\bf B}(p_i)\,$ r.v.s one has $\,X_i^*\equiv 0$, and (\ref{314}) coincides with (\ref{be}). 

	In the case of i.i.d.r.v.s (\ref{314}) becomes
 $$ \d(S_n;\pi_\l) \le (1\!-\!e^{-\l}) \E|X\!-\!X^*| \,. $$
  Here $\,X^*\,$ may be chosen independent of $\,X$, although one would prefer to define $\,X\,$ and $\,X^*\,$ on a common probability space in order to make $\,\E |X\!-\!X^*|\,$ smaller. 

	A generalisation of (\ref{RN}) to the case of independent integer-valued r.v.s has been given by Novak \cite{N18}.\\

   \noindent{\bf Example 2.2}. Let $\,\xi,X_1,X_2,...\,$ be i.i.d.r.v.s with geometric $\,{\bf \Gamma}_0(p)\,$ distribution: 
	$$ \p(\xi\!=\!m) = (1\!-\!p)p^m\quad (m\!\ge\!0). $$ 
	Then $\,S_n\,$ is a negative Binomial $\,{\bf N\!B}(n,p)\,$ r.v.. 

  Set $\,r\!=\!p/(1\!-\!p).$ Vervaat \cite{V69} has shown that 
$\,\d(S_n;\pi_\l)\le r,$ while Romanowska \cite{Ro77} has noticed that $\,\d(S_n;\pi_\l)\le r/\sqrt{2}.$ 
		Roos \cite{R03JSPI} has shown that  
  \b \label{R03}  \d({\bf N\!B}(n,p);{\bf\Pi}(np)) &\le& \min\{3r\!/4e;nr^2\},\\ 
		 \label{R03G}	d_{_G}({\bf N\!B}(n,p);{\bf\Pi}(np)) &\le& nr^2. \e

  It is easy to see that $\,\p(X_i^*\!=\!m)=(m\!+\!1)p^m(1\!-\!p)^2\,.$ Hence
 \beq                                                       \label{319}
 X_i^*\st X_i+\xi, 
 \eeq 
and $\,\E |X\!-\!X^*| = p/(1\!-\!p)$. Note that 
 $$ \l=n\E\xi = nr,\ \ d_{_G}(X;X^*) = \E\xi=r. $$ 
Theorem \ref{P38} entails 
 \b \label{317}
 \d(S_n;\pi_\l) &\le& (1\!-\!e^{-nr})r,\\ \label{2G} d_{_G}(S_n;\pi_\l) 
 &\le& \min\!\Big\{1;\frac{_4}{^3}\sqrt{2/enp}\,\Big\} \,nr^2. 
 \e 
   Inequality (\ref{317}) has been established in \cite{B87}, p. 758; estimate (\ref{2G}) is from \cite{N11}, formula (4.53). \hspace*{\fill} $\Box$\\

      {\bf Shifted Poisson approximation}. A number of authors dealt with shifted Poisson approximation to the distribution of a sum $\,S_n\,$ of integer-valued r.v.s (see \cite{BC02,N18} and references therein). 
   Let 
	$$ \l=\E S_n,\ \sigma^2 = \hbox{var}\, S_n,\ a = [\l\!-\!\sigma^2],\ 
     b = \{\l\!-\!\sigma^2\},\ \mu = \sigma^{2}\!+b, $$ 
where $\,[x]\,$ and $\,\{x\} = x\!-\![x]\,$ denote the integer and the fractional parts of $\,x$. 

	Barbour \& $\rm\check{C}$ekanavi\v{c}ius \cite{BC02} have shown that 
 \beq 				\label{P14} 
 \d(S_n;a\!+\!\pi_{\mu}) \le (1\wedge\sigma^{-2}) 
 \Big(b+d_n\sum^n_{i=1} \psi_i\Big) + \p(S_n\!<\!a), 
 \eeq 
where $\,d_n = \max_{i\le n}\d(S_{n,i};S_{n,i}\!+\!1),$ $\,S_{n,i} = S_n\!-\!X_i,$ $\,\psi_i = \sigma^2_i \E X_i(X_i\!-\!1) + |\E X_i\!-\!\sigma^2_i|\E (X_i\!-\!1)(X_i\!-\!2) + \E |X_i(X_i\!-\!1)(X_i\!-\!2)|,$ $\,\sigma^2_i=\hbox{var} X_i$. 

   In the Binomial case (i.e., $\,\L(S_n)={\bf B}(n,p)$) the r.-h.s. of (\ref{P14}) is $\,O\big(\sqrt{p/n}\,+1/np\big)$. 	Further reading on the topic is \cite{N18}. 

	An estimate of the accuracy of shifted Poisson approximation to the distribution of a random sum of i.i.d. integer-valued r.v.s has been presented by R\"ollin \cite{Ro05}.

					\subsection{Dependent integer-valued r.v.s}\label{2.4}

  Let $\,X_1,...,X_n\,$ be (possibly dependent) non-negative integer-valued r.v.s. Set $\,p_i = \p(X_i\!=\!1|X_1,...,X_{i-1}).$ A generalisation of (\ref{LeCam}$^*$), (\ref{S75}) has been given by Serfling \cite{S75}: 
 $$  
 \d(S_n;\pi_\l) \le \sum^n_{i=1} \(\E^2p_i + \E|p_i\!-\!\E p_i| + \p(X_i\!\ge\!2)\), 
 \eqno(\ref{LeCam}^+) $$ $$ \ \ 
 d_K(S_n;\pi_\l) \le \sum^n_{i=1} \( \frac{\,_2}{\,^\pi} \E^2p_i + \E|p_i\!-\!\E p_i| + \p(X_i\!\ge\!2) \). \eqno(\ref{S75}^+) $$ 

	Below we present a generalisation of Theorem \ref{T37}. 

   Let $\,\{X_a, a\!\in\!J\}\,$ be a family of r.v.s taking values in $\,\Z$. 
Suppose one can choose the ``neighborhoods" $\,\{B_a\}\,$ so that r.v.s $\,\{X_b, b\!\in\!J\!\setminus\! B_a\}\,$ are independent of $\,X_a\,.$ 
   We call this assumption the ``local dependence''  condition.

   Let $\,\L(\pi_\l)\,$ denote a Poisson $\,{\bf\Pi}(\l)\,$ r.v.. Set 
 \begin{eqnarray*}  
 \delta_1^* &=& \sum_{a\in J} \sum_{b\in B_a\backslash \{a\}} \!\E X_a \E X_b
 \,,\ \ \delta_4 = \sum_{a\in J} d_{_G}(X_a;X_a^*)\E X_a ,  
 \end{eqnarray*} 
and let $\,\delta_1, \delta_2, \delta_3\,$ be defined as in Theorem \ref{T37}. 
  Theorems \ref{T39} and \ref{T310} are from \cite{N11}, ch. 4. 

   \begin{theorem} \label{T39} 
If $\,\{X_b, b\!\in\!J\!\setminus\!B_a\}\,$ are independent of $\,X_a,$ then
 \beq                                               		\label{315}
 \d(S_n;\pi_\l) \le \frac{1\!-\!e^{-\l}}{\l} \(\delta_1^* + \delta_2 + \delta_4\). \eeq \end{theorem}

   In Theorem \ref{T310} we drop the local dependence condition assumed in Theorem \ref{T39}.  

   \begin{theorem} \label{T310}  
Denote $\,\delta_5 = \sum_{a\in J} \E X_a(X_a\!-\!1) \1\{X_a\!\ge\!2\}.$ Then
 \beq                                               		\label{316}
 \d(S_n;\pi_\l) \le \frac{1\!-\!e^{-\l}}{\l} \Big(\delta_1 + \delta_2 + \delta_5\Big) 
 + \min\{1;\sqrt{2/e\l}\,\} \delta_3 . \eeq    \end{theorem} 

  Ruzankin \cite{R10} presents an estimate of the accuracy of Poisson approximation to $\,\E h(S_n),$ where $\,h\,$ is an unbounded function.\\ 

  \noindent{\it Open problem}.\\ 
  \noindent2.4. Improve the constants in (\ref{315}), (\ref{316}).

							\subsection{Asymptotic expansions}\label{2.6}

	Let $\,X_1,...,X_n\,$ be independent Bernoulli $\,{\bf B}(p_i)\,$ r.v.s, and let $\,\pi_\l\,$ be a Poisson random variable. 
	
	Formal expansions of $\,\p(S_n\!\le\!x)\,$ have been given by Uspensky \cite{U31}, see also Franken \cite{Fr}. 
	Herrmann \cite{H65}, Shorgin\index{Shorgin} \cite{Sho77} and Barbour \cite{B87} 
present full asymptotic expansions with explicit estimates of the error terms. 
  Kerstan \cite{K64}, Kruopis \cite{K86} and $\rm\check{C}$ekanavi\v{c}ius \& Kruopis \cite{CK00} present first-order asymptotic expansions. 
   Asymptotic expansions for $\,\E h(S_n)-\E h(\pi_\l)\,$ in the case of independent 0-$\!$1 r.v.s $\,\{X_k\}\,$ and unbounded function $\,h\,$ have been given by Barbour et al. \cite{BCC95} and Borisov \& Ruzankin \cite{BR02}. 

  The formulation of the full asymptotic expansions is cumbersome and will be omitted. We present below first-order asymptotics	of $\,\E h(S_n)\,$ for particular classes of functions $\,h$. 

   Of special interest are indicator functions $\,h(\cdot)=\1\{\cdot\!\in\!A\},\ A\!\subset\!\Z$. Denote\index{Asymptotic expansions!in the Poisson limit theorem} 
 \bb 
 Q_\l(A) \!&=&\! \Big[ \p(\pi_\l\!\in\!A) + \p(\pi_\l\!+\!2\!\in\!A) -
           2\p(\pi_\l\!+\!1\!\in\!A) \Big]\Big/2,\\ 
 \ve \!&=&\! \min\!\left\{1; \(2\pi[\l\!-\!p^*_n]\)^{-1/2} + 
 2\delta/(1\!-\!p^*_n/\l)\right\},\ \ p^*_n=\max_{i\le n}p_i. \ee 

	Let $\,\pi_\l^\star\,$ denote a random variable with distribution (\ref{star}). Then 
	$$\,Q_\l(A) = [\p(\pi_\l^\star\!\in\!A)-\p(\pi_\l\!+\!1\in\!A)]/2\l\,$$ 
(see \cite{N11}, ch. 4). 

  The following result from \cite{N11}, ch. 4, sharpens (13) in \cite{H65} and the bound of Corollary 2.4 in \cite{B87} (Corollary 9.A.1 in \cite{BHJ}).

  \begin{theorem} \label{asy} 
	Let $\,X_1,...,X_n\,$ be independent Bernoulli r.v.s, $\,\L(X_i)={\bf B}(p_i)$. Then 
 \beq                                              \label{asym} 
 \left|\p(S_n\!\in\!A) - \p(\pi_\l\!\in\!A) + Q_\l(A) \sum_{i=1}^n p_i^2 \right|
 \le 2\delta^*\ve + 2\delta^2\,, 
 \eeq 
where $\,\delta = \l^{-1}(1\!-\!e^{-\l}) \sum^n_{i=1} p_i^2\,,$ 
$\,\delta^* = \l^{-1}(1\!-\!e^{-\l}) \sum^n_{i=1} p_i^3$. \end{theorem}

  Recall that $\,\l_k=\sum_{i=1}^n p_i^k\ \ (k\!\ge\!2).$ Denote 
$$ \Delta h(\cdot) = h(\cdot+1)-h(\cdot). $$

	\begin{theorem} \label{BoRu} {\rm\cite{BR02}} 
If $\,\E |h(\pi_\l)|\pi_\l^4 <\infty,$ then 
	\b \nonumber 
	&& \Big| \E h(S_n) - \E h(\pi_\l) + \l_2 \E\Delta^2 h(\pi_\l)/2	\Big|\\ 
	&& \le \frac{e^{p^*_n}}{(1\!-\!p^*_n)^2} \Big( \l_3 \E|\Delta^3 
     h(\pi_\l)|/3 + \l_2^2 \E|\Delta^4 h(\pi_\l)|/8 \Big). \label{BoRuz} 
	\e \end{theorem}

  Note that the assumption $\,\E|\Delta^k h(\pi_\l)|\!<\!\infty\,$ is equivalent to $\,\E\pi_\l^k|h(\pi_\l)|\!<\!\infty\ (k\!\in\!\N),$ see Proposition 1 in \cite{BR02}. 	Borisov \& Ruzankin (\cite{BR02}, Lemma 2) have showed also that 
  $$ \sup_k \p(S_n\!=\!k)/\p(\pi_\l\!=\!k) \le (1\!-\!p_n^*)^2 . $$ 

  Asymptotic expansions for $\,\E h(S_n)-\E h(\pi_\l),$ where $\,\{X_k\}\,$ are non-negative integer-valued random variables and function $\,h\,$ is either bounded or grows at a  polynomial rate, are presented in Barbour \cite{B87}. 
  Asymptotic expansions for $\,\E h(S_n)-\E h(\pi_\l),$ where $\,\|h\|_1=1$, have been given by Barbour \& Jensen \cite{BJ89}.\\

      {\bf Unit measure (signed measure) approximations}. 
	 A number of authors evaluated the accuracy of unit measure (signed measure) approximation to the distribution of a sum $\,S_n\,$ of independent Bernoulli r.v.s (see, e.g., \cite{B88,BX99,BC02}). 
  In particular, Borovkov \cite{B88} has generalised inequality (\ref{LeCam}). 
	Note that asymptotic expansion (\ref{asym}) is an example of a unit measure approximation. 
	
	Denote by $\,P_n\,$ the distribution corresponding (with some abuse of notation) 
to $\,\pi_{\l+\l_2}+2\pi_{-\l_2/2}\,$ (i.e., $\,P_n\,$ is a convolution of $\,{\bf\Pi}(\l\!+\!\l_2)\,$ and a Poisson unit measure with parameter $\,-\l_2/2\,$ on $\,2\Z$). In the assumption that $\,\theta\!<\!1/2\,$	Barbour \& Xia (\cite{BX99}, Theorem 4.1) have shown that 
	$$ \d(\L(S_n);P_n) \le \l_3/\l(1\!-\!2\l_2)\sqrt{\l\!-\!\l_2\!-\!p^*_n}\,. $$ 
	
	$\rm\check{C}$ekanavi\v{c}ius \& Kruopis \cite{CK00} present an estimate of the accuracy of unit measure approximation in terms of the Gini-Kantorovich distance: 
	$$ 
	d_{_{G}}(\L(S_n);Q_n) \le C\l_*^{-1}\l_2(1+(\l/\l_2)^2), 
	$$ 
where $\,C\,$ is an absolute constant, $\,\l_* = \max\{1;\l\!-\!\l_2\}\,$ and $\,Q_n\,$ (with some abuse of notation) corresponds to $\,\pi_{\l-\l_2/2}-\pi_{-\l_2/2}\,$ ($\,Q_n\,$ is a convolution of $\,{\bf\Pi}(\l\!-\!\l_2/2)\,$ and a Poisson unit measure with parameter $\,-\l_2/2\,$ on $\,-\Z$). Note that 
  $\,\sum_k kQ_n(k)=\E S_n,$ $\,\sum_k(k\!-\!\l)^2P_n(k)=\hbox{\rm var\,}S_n.$

	Barbour \& $\rm\check{C}$ekanavi\v{c}ius \cite{BC02} present a unit measure approximation to the distribution of a sum of independent integer-valued r.v.s.

		\subsection{Sum of a random number of random variables}\label{2.5} 

  Let $\,\nu,X,X_1,X_2,...\,$ be independent non-negative random variables, where r.v. $\,\nu\,$ takes values in $\,\Z,$ $\,X,X_1,X_2,...\,$ are i.i.d. random variables. 
	
	Set $$\,S_\nu = X_1+...+X_\nu\,.$$ A natural task is to evaluate the accuracy of Poisson approximation to $\,\L(S_\nu).$ 

   We consider first the case where $\,X,X_1,X_2,...\,$ are Bernoulli $\,{\bf B}(p)\,$ r.v.s.\\ 
 
  Denote $\,\bar\nu := \E\nu$. Then $\,\E S_\nu = p\bar\nu.$ 
	
	Let $\,{\cal F}_2\,$ denote the class of functions $\,h:\Z\!\to\!\R\,$ such that 
$\,||\Delta^2h||\!\le\!1$, and set 
  $$ d_2(X;Y) = \sup_{h\in{\cal F}_2} |\E h(X)-\E h(Y)|. $$ 
Logunov \cite{Log90} points out that $\,\d(X;Y) \le d_2(X;Y),$ and shows that 
  $$ d_2(S_\nu;\pi_{p\bar\nu}) \le p^2d_*(\nu;\pi_{p\bar\nu}), $$ 
where $\,d_*(X;Y) = \sum_{k\ge1}k(k\!-\!1)|\p(X\!=\!k)-\p(Y\!=\!k)|/2.$ 
  Note that $\,d_*(X;Y)\ge d_{_G}(X;Y).$ 

   Yannaros \cite{Y91} has shown that 
	\beq                                                           \label{Y}
  \d(\pi_\l;\pi_\mu) \le \min\{|\sqrt{\l}-\sqrt{\mu}|;|\l\!-\!\mu|\}. 
  \eeq 
  The first term in (\ref{Y}) has been improved by Roos \cite{R03JSPI}: 
$$ \d(\pi_\l;\pi_\mu) \le \sqrt{\frac{_2}{^e}}\,\Big|\sqrt{\l}-\sqrt{\mu}\,\Big|. 
	   \eqno(\ref{Y}^*) $$ 
	Note that the second term in the r.-h.s. of (\ref{Y}) is a consequence of the trivial inequality 
	$$ \d(\pi_\l;\pi_\mu) \le 1-\exp(-|\l\!-\!\mu|) \eqno(\ref{Y}^\star) $$ 
that follows by defining $\,\pi_\l\,$ and $\,\pi_\mu\,$ on a common probability space (cf. (4.10) in \cite{N11}). 

	It is easy to see that 
  \bb 
  \d(S_\nu;\pi_\l) &\le& \sum_k\nolimits \p(\nu\!=\!k) \d(S_k;\pi_\l),\\ 
  \d(S_k;\pi_\l) &\le& \d(S_k;\pi_{kp}) + \d(\pi_{kp};\pi_\l).\nonumber \ee 
Using these inequalities and (\ref{Y}), Yannaros \cite{Y91} has shown that 
 \beq                                         \label{Yan2}
 \d(S_\nu;\pi_{p\bar\nu}) \le \min\!\Big\{ \frac{p}{2\sqrt{1\!-\!p}}\, ; 
 (1\!-\!\E e^{-p\nu})p \Big\} + \min\Big\{ p\E|\nu\!-\!\bar\nu| ; 
 \sqrt{p\frac{\hbox{\rm var\,}\nu}{\bar\nu}}\,\Big\}. \eeq 
   The term $\,\min\!\big\{ p/2\sqrt{1\!-\!p}\, ; (1\!-\!\E e^{-p\nu})p \big\}\,$ 
in (\ref{Yan2}) is inherited from (\ref{Rom}) and (\ref{be}). 

	The right-hand side of (\ref{Yan2}) can be sharpened using (\ref{RN}), (\ref{Y}$^\star$) and (\ref{Y}$^*$): 
  \beq \label{YN} 
	\d(S_\nu;\pi_{p\bar\nu}) \le 3p/4e + 2(\delta^*\!+\!\delta^2) + 
	\min\Big\{ \Big(1\!-\!\E e^{-p|\nu-\bar\nu|}\Big) ; 
	\E|\sqrt{\nu}-\sqrt{\bar\nu}\,|\sqrt{2p/e}\, \Big\}. \eeq 
Note that $\,\E|\sqrt{\nu}-\sqrt{\bar\nu}\,| \le \min\{ \bar\nu^{-1/2} \sqrt{\hbox{\rm var\,}\nu}\,; \bar\nu^{-3/2}\hbox{\rm var\,}\nu \}.$\\ 
	
  {\bf Mixed Poisson distribution}. 
	A number of authors (see, e.g., Roos \cite{R03JSPI}) have evaluated the accuracy of Poisson approximation to the {\it mixed Poisson distribution}, i.e., the distribution of the r.v. $\,\pi_\nu,$ where $\,\L(\pi_t)={\bf\Pi}(t)$, r.v. $\,\nu\,$ takes values in $\,[0;\infty)\!:$ 
  $$ \p(\pi_\nu\!=\!m) = \int_0^\infty \!\p(\pi_y\!=\!m)\p(\nu\!\in\!dy) 
	\qquad(m\!\ge\!0). $$ 
  If $\,\{X_i\}\,$ are Poisson $\,{\bf\Pi}(\l)\,$ r.v.s, then $\,S_\nu\st\pi_{\l\nu}\,$ is a mixed Poisson random variable. 

	Denote by $\,{\bf N\!B}(n,p)\,$ the negative Binomial distribution: 
$\,\L(S_n) = {\bf N\!B}(n,p)\,$ if 
	$$ \p(S_n\!=\!i) = {i\!+\!n\!-\!1\choose i}(1\!-\!p)^np^i\ \qquad (i\!\ge\!0). $$ 
  The negative Binomial distribution $\,{\bf N\!B}(t,p)\,$ is a  mixed Poisson distribution with 
  $$ \p(\nu\!\in\!dy)/dy = r^ty^{t-1}e^{-yr}/\Gamma(t) \qquad(y\!>\!0), $$ 
where $\,r\!=\!p/(1\!-\!p),$ $\,\Gamma(y) = \int_0^\infty x^{y-1}e^{-x}dx.$ 

  Roos \cite{R03JSPI} presents estimates of the accuracy of Poisson approximation to the mixed Poisson distribution with a correct constant at the leading term. 

  $\,$ 

		{\bf Sum of 0-$\!$1 random variables till the stopping time}.	
		We now consider the situation	where r.v. $\,\nu\,$ depends on $\,\{X_i\}$. 

		Let $\,X,X_1,X_2,...\,$ be i.i.d. non-negative integer-valued r.v.s.  
		Set $\,S_0\!=\!0,$ $$\,S_n=X_1+...+X_n\quad (n\!\ge\!1),$$ 
and let $\,\mu(t)\,$ denote the stopping time: 
$$ \mu(t) = \max\{n\!\ge\! 0: S_n\!\le\!t\}. $$ 
	  Theorems \ref{7T4}--\ref{7T5} below are cited from see \cite{N11}, ch. 3. They provide estimates of the accuracy of Poisson approximation to the distribution of the number 
 \beq                                                           \label{7b10}
 N_t(x) = \sum_{j=1}^{\mu(t)} \1\{X_j\!\ge\! x\} + \1\{t\!-\!S_{\mu(t)}\!\ge\!x\} 
 \eeq
of exceedances of a ``high'' level $\,x\!\in\![0;t]\,$ till $\,\mu(t)$. 

   Note that 
 $$ \{N_t(x)\!=\!0\} = \{M_t\!<\!x\}, $$ 
where 
 \beq                                               \label{7b1}
 M_t = \max\{t\!-\!S_{\mu(t)} ; \max_{1\le i\le \mu(t)} X_i\}
 \eeq 
is the largest observation among $\,\{X_1,...,X_{\mu(t)}, t-S_{\mu(t)}\}$. 

  Let $\,X_{k,t}\,$ denote the $k^{\rm th}$ largest element among
$\,\{X_1,...,X_{\mu(t)}, t-S_{\mu(t)}\}$. Then $$\,\{X_{k,t}\!<\!x\}=\{N_t(x)\!<\!k\}.$$ 

  The topic has applications in finance. For instance, suppose a bank has opened a credit line for a series of operations, and the total amount of credit is $\,t\,$ units of money. The cost of the $i$-th operation is denoted by $\,X_i\,.$ 
  What is the probability that the bank will ever pay $\,x\,$ or more units of money at once? that there will be a certain number of such payments? 
  Information on the asymptotic properties of the distribution of random variables $\,M_t\,$ and $\,N_t(x)\,$ can help to answer these questions. 

   Let $\,\{X_i^\<,i\!\ge\!1\},\,\{X_j^\>,j\!\ge\!1\}\,$ be independent r.v.s
with the distributions 
   $$ \L(X^\<) = \L(X|X\!<\!x),\ \L(X^\>) = \L(X|X\!\ge\!x). $$ 
  We set $\,p_x = \p(X\!\ge\!x),$  
 \bb 
 S_0(k)=0,\ \ 
 S_m(k) \!&=&\! \sum^k_{i=0} X_i^\> + \sum^m_{i=k+1} X_i^\< \ \ \ (m\!\ge\!1).
 \ee
   Let $K_*,\,K^*$ denote the end-points of $\,\L(X),$ and set 
 \bb
 \tau_k &=& \tau_k'\!-\!k\,,\ \tau_k' = \min\{n: S_n(k)>t\!-\!x\},\\ 
 \l_k &\equiv& \l_k(t,x,k) = p_x (t\!-\!x\!-\!k\E X^\>)/\E X^\<.
 \ee  

In Theorems \ref{7T4}--\ref{7T5} we assume the following condition:\\ {\it there exist constants $\,D\!<\!\infty\,$ and $\,D_*\!\in\!\(K_*;K^*\)\,$ such that} 
 \beq                                                       \label{7b9}
 \int_x^\infty\! \p(X\!\ge\!y)dy \le D\p(X\!\ge\!x)\qquad(x\!\ge\!D_*). \eeq 

   Condition (\ref{7b9}) means the tail of $\,\L(X)\,$ is light (cf. (3.15) in \cite{N11}). Inequality (\ref{7b9}) holds if function $\,g(x) = e^{cx}\p(X\!\ge\!x)\,$ is not increasing as $\,x\!>\!1/c\ \,(\exists c\!>\!0)$.
   The equality in (\ref{7b9}) for all $\,x\!\ge\!0\,$ may be attained only if $\,\L(X)\,$ is exponential with $\,\E X= D$.

   \begin{theorem} \label{7T4}  
For any $\,k\!\in\!\Z,$ as $\,t\to\infty,$ 
 $$
 \sup_{x\in B_+(t)}\! \left| \p(N_t(x)\!=\!k) - \p(\pi_{\l_k}\!\!=\!k) -
 \sum_{r=0}^{k-1}\! \(\p(\pi_{\l_k}\!\!=\!r)-\p(\pi_{\l_{k-1}}\!\!=\!r) \) \right|
 \!=\! O(1/t), $$
where $\,B_+(t)=(K_*;K^*\wedge t/(k\!+\!2)).$ \end{theorem}   

    Let $\,\pi(t,x)\,$ denote a Poisson r.v. with parameter $\,p_xt/\E X$.

   \begin{theorem} \label{7T5}  
For any $\,k\!\in\!\Z,$ as $\,t\to\infty,$
 $$ \sup_{K_*<x<K^*}\! \left| \p(N_t(x)\!=\!k) - \p(\pi(t,x)\!\!=\!k) \right|
 = O(t^{-1}\ln t). $$ \end{theorem}

	One can show that $\,N_t(x)\,$ is ``small'' when $\,x\,$ is ``large'': 
 $$ \sup_{x\ge \sqrt{t}} \p(N_t(x)\!\ge\!1) \le 
 q^{\sqrt{t}} \qquad(\exists q\!\in\!(0;1)). $$ 
  Theorem 3.7 in \cite{N11} presents asymptotic expansions for $\,\p(N_x(t)\!=\!k)$. The asymptotic expansions for $\,\L(M_t)\,$ are available under a weaker moment assumption (cf. \cite{N11}, ch. 3). 

   $\,$  

   \noindent{\bf The number of intervals between consecutive jumps of a Poisson process.} 
Consider a Poisson jump process $\,\{\pi_\l(s), s\!\ge\!0\}\,$ with parameter $\,\l\!>\!0,$ and let $\,\eta_i\,$ denote the moment of its $i^{\rm th}$ jump. 
	Set $\,X_i=\eta_i-\eta_{i-1}\,.$ Then $\,N_t(x)\,$ is the number of intervals between consecutive jumps with lengths greater or equal to $\,x$. If the points of jumps represent catastrophic/rare events, then $\,N_t(x)\,$ can be interpreted as the number of ``long'' intervals without catastrophes.

    Let $\,\pi_{t,x}\,$ be a Poisson r.v. with parameter $\,t\l e^{-\l x}\,.$
Then for any $\,k\!\in\!\Z,$ as $\,t\to\infty,$ 
 \beq                                           \label{7}
 \sup_{0<x<t} \left| \p(N_t(x)=k) - \p(\pi_{t,x}=k) \right|
 = O\big(t^{-1}\ln t\big) 
 \eeq (cf. (3.12) in \cite{N11}).\\

 \noindent{\it Open problems}.\\ \index{open problems}
 \noindent2.5. Will asymptotic expansions for $\,\L(N_x(t))\,$ hold under a weaker moment assumption?\\ 
 \noindent2.6. Generalise  the results of Theorems \ref{7T4}--\ref{7T5} to the case of
 $$ N_t(x) = \sum_{j=1}^{\mu(t)} Y_i\1\{X_j\!\ge\! x\} +
 Y_{\mu(t)+1} \1\{t\!-\!S_{\mu(t)}\!\ge\!x\}, $$
where $\,\{(X_i,Y_i)_{i\ge1}\}\,$ is a sequence of i.i.d. pairs of r.v.s, $\,Y_i\!>\!0$. 

							\section{Applications}\label{3}  

	Applications of the theory of Poisson approximation to meteorology, reliability theory and extreme value theory have been discussed in \cite{BalK,H67,LLR,N11}. In this section we present a number of results that are not fully covered in existing surveys. 
	
					\subsection{Long head runs}\label{3.1}

	 Let $\,\{\xi_i,i\!\ge\!1\}\,$ be a sequence of 0-$\!$1 random variables. 
	
   We say a head run (a series of 1's) starts at $\,i=1\,$ if $\,\xi_1=1$; a series starts at $\,i>1\,$ if $\,\xi_{i-1}=0,\, \xi_i=1$.
   If $\,\xi_{i-1}\!=\!0,\, \xi_i=\!...\!=\xi_{i+k-1}=1,$ we say the head run is of length $\,\!\ge\! k$. 
	
	For instance, if $\,n=5\,$ and $\,\xi_1\!=\!\xi_2\!=\!\xi_3\!=\!1, \xi_4\!=\!0, \xi_5\!=\!1,$ there is one series (head run) of length 3 and one series of length 1. 

   Denote $$\,A_0 = \{\xi_{1}=\!...\!=\xi_{k}=1\},\ 
\,A_i = \{\xi_i\!=\!0,\xi_{i+1}=\!...\!=\xi_{i+k}=1\}\quad (i\!>\!1).$$ 
   Then 
 $$ W_n(k) = \sum_{i=0}^{n-k} \1\{A_i\} \qquad(n\!\ge\!k\!\ge\!1) $$
is the number of head runs of length $\,\!\ge\! k\,$ among $\,\xi_1,...,\xi_n\,$ (NLHR).  

   Set  
 \beq                                                   \label{7b0}
 L_n = \max\{k: \xi_{i+1}=...=\xi_{i+k}=1\ (\exists i\!\le\!n\!-\!k)\}. 
 \eeq 
$\,L_n\,$ is the {\it length of the longest head run} (LLHR)\index{LLHR} among 
$\,X_1,...,X_n$. Obviously, $$\,\{L_n\!<\!k\} = \{W_n(k)\!=\!0\}.$$ 

   The problem of approximating the distribution of LLHR is a topic of active research; it has applications in reliability theory and psychology (cf. \cite{BalK,N11}).\\ 

   Let $\,\{\xi_i,i\!\ge\!1\}\,$ be i.i.d. Bernoulli $\,{\bf B}(p)\,$ r.v.s,
$\,p\!\in\!(0;1),$ and let $\,\pi_\l\,$ denote the Poisson $\,{\bf\Pi}(\l)\,$ r.v.. 
   Theorem \ref{T37} with $\,B_i = [i\!-\!k;i\!+\!k]\,$ and 
$$\,\l \equiv \l(n,k,p) = p^k(1\!+\!(n\!-\!k)(1\!-\!p)) $$ yields the following

 \begin{corollary} \label{Cb30} As $\,n\!\ge\!k\!\ge\!1,$ 
 \beq                                           \label{7b30}
 \d(W_n(k);\pi_\l) \le (1\!-\!e^{-\l}) (2k\!+\!1) p^k . \eeq
 \end{corollary} 

   An open question is if estimate (\ref{7b30}) can be improved. Note, for instance, that (\ref{7b30}) does not yield (\ref{7c}) even for $\,j\!=\!0$.\\ 

   There is a close relation between $\,N_t(x)\,$ and $\,W_n(k)$. Let $\,\eta_0=0$,
 $$
 \eta_i = \min\{k\!>\!\eta_{i-1}\!: \xi_k=0\},\ X_i=\eta_i-\eta_{i-1}\ \ (i\!\ge\!1).
 $$
Then
 \beq                                               \label{7b}
 W_n(k) = \sum_{j=1}^{\mu(t)} \1\{X_j\!-\!1\!\ge\!k\} 
 + \1\{n\!-\!\eta_{\mu(n)}\!\ge\!k\}. 
 \eeq
Hence $$ W_{n-1}(k)=N_n(k\!+\!1). $$  

    Denote $\,\l_k=n(1\!-\!p)p^k\,$. Theorem \ref{7T5} entails

   \begin{corollary} \label{7C6} For any $\,j\in\Z,$ as $\,n\to\infty,$
 \beq 																\label{7c}
 \max_{1\le k\le n} \left| \p(W_n(k)=j) - \p(\pi_{\l_k}=j) \right|
 = O\(n^{-1}\ln n\). \eeq \end{corollary}  
 
   According to Theorem 3.13 in \cite{N11}, the rate $\,n^{-1}\ln n\,$ in (\ref{7c}) cannot be improved.\\

				{\bf The number of long non-decreasing runs}. 
	Let $\,\xi_i = \1\{Y_i\!\le\!Y_{i+1}\},$ where $\,\{Y_i\}\,$ are i.i.d.r.v.s with a continuous d.f.. Then NLHR $\,W_n(k)\,$ is the number of non-decreasing runs of length $\,\!\ge\!k\,$ (NLNR), and LLHR is the length of the longest non-decreasing run (LLNR) among $\,Y_1,...,Y_{n+1}$. 
	 We denote LLNR by $\,L_n^+\,$ and NLNR by $\,W_n^+(k)$. 

   The topic concerning LLNR and NLNR has applications in finance. 
   It is well known that prices of shares and financial indexes evolve in cycles of growth and decline. 
	Knowing the asymptotics of $\,L_n^+\,$ and $\,W_n^+(k)\,$ can help evaluating the length of the longest period of continuous growth/decline of a particular financial instrument as well as the distribution of the number of such long periods. 

   Pittel \cite{P81} has proved a Poisson limit theorem for NLNR (see also Chryssaphinou et al. \cite{CPV} concerning the case of a Markov chain). 

   We proceed with the case of i.i.d.r.v.s with a continuous d.f.. 
   Note that $\,\L(\xi_i)={\bf B}(1/2)\,$ and $\,\p(Y_1\!\le...\le\!Y_{k+1}) = 1/(k\!+\!1)!$. Set $\,\l_{n,k} =\E W_n^+(k).$ Then 
	$$ \l_{n,k} = 1/(k\!+\!1)! + (n\!-\!k)/k!(k\!+\!2). $$ 
   Theorem \ref{T37} with $\,B_i = [i\!-\!k\!-\!1;i\!+\!k\!+\!1]\,$ yields the following

 \begin{corollary} \label{Cinc} As $\,n\!\ge\!k\!\ge\!1,$ 
$$ \d(W_n^+(k);\pi_{\l_{n,k}}) \le (1\!-\!e^{-\l_{n,k}}) (2k\!+\!3)/(k\!+\!1)!. $$  \end{corollary}

	The accuracy of compound Poisson approximation to the distribution of the number of non-decreasing runs of fixed length has been evaluated by Barbour \& Chryssaphinou \cite{BC01}, p. 982 (continuous d.f.) and Minakov \cite{Min15} (discrete d.f.). 
	Concerning the asymptotics of LLNR, see \cite{CF96,N11} and references therein.\\ 

 \noindent{\it Open problem}.\\ \index{open problems} 
 \noindent{3.1}. Improve the estimates of Corollary \ref{Cb30} and Corollary \ref{Cinc}.\\ 
 \noindent{3.1}. Derive (\ref{PEx3})-type (i.e., uniform in $\,k$) estimates of the accuracy of (possibly shifted) Poisson approximation to $\,\L(W_n(k))\,$ and $\,\L(W_n^+(k))$.

					\subsection{Long match patterns}\label{3.2}

   Closely related to the number of long head runs is the number of long match patterns (NLMP) between sequences of independent r.v.s. 
   Information on the distribution of NLMP and the length of the longest match pattern (LLMP) can help recognising ``valuable'' fragments of DNA sequences (see \cite{AGG,AGW90,Mi08,M90,Ne94}). 

  In this section we present results on the accuracy of Poisson approximation to the distribution of NLMP. 
	Theorems \ref{LMP}, \ref{7bT14} and Lemma \ref{LMP-L1} below have been established by the author (see \cite{N11}, ch. 4). 

   Let $\,X,X_1,..., X_m\,,\,Y,Y_1,..., Y_n\,$ be independent non-degenerate random variables taking values in a discrete state space $\,\bf A$. Denote $\,(k\!\in\!\N)$ 
 \bb 
 T_{ij} \!&=&\! \1\{X_{i+1}\!=\!Y_{j+1},\ldots, X_{i+k}\!=\!Y_{j+k}\},\\
 \tilde T_{ij} \!&=&\! T_{ij}(k)\1\{X_i \!\ne\! Y_j\},\\ 
 T_{ij}^* \!&=&\! \tilde T_{ij}\ \ (i\!\ge\! 1,j\!\ge\!1),\ \ 
 T_{ij}^* = T_{ij}\ \ (i\!=\!0\ \hbox{or}\ j\!=\!0). 
 \ee 
    Then  
 $$
 M_{m,n}^* = \max\Big\{k\!\le\!\min(m,n): \max_{(i,j)\!\in\!J} T_{ij}=1\Big\}
 $$
is the {\it length of the longest match pattern} between 
$\,(X_1\ldots X_m)\,$ and $\,(Y_1\ldots Y_n)$.  

   LLMP $\,M_{m,n}^*\,$ is a 2--dimensional analog of LLHR $\,L_n\,.$ 
If $\,{\bf A}=\{0,1\}\,$ and $\,Y_1=...=Y_n=1,$ then $\,M_{n,n}^* = L_n.$ 

   Given $\,m\!\ge\!k,\,n\!\ge\!k$, let 
 $$
 J \equiv J(k,m,n) = \{(i,j): 0\!\le\!i\!\le\!m\!-\!k, 
 0\!\le\!j\!\le\!n\!-\!k\}. $$ 
   Denote by
 $$ W_{m,n} \equiv W_{m,n}(k) = \sum_{(i,j)\in J}T_{ij}^* $$
the {\it number of long match patterns} (patterns of length $\!\ge\! k$). Then 
 $$ \{M_{m,n}^*\!<\!k\} = \{W_{m,n}\!=\!0\}. $$ 

    In the rest of this section we assume that r.v.s $\,X,X_1,..., X_m\,,\,Y,Y_1,..., Y_n\,$ are identically distributed. 
  We set 
	$$ \l \equiv \l_{k,m,n} = \E W_{m,n},\ \ m'=m\!-\!k+\!1,\ \ n'=n\!-\!k\!+\!1. $$ 
Then $\,\l = (m'\!-\!1)(n'\!-\!1)(1\!-\!p)p^k + (m'\!+\!n'\!-\!1)p^k.$ 

Denote 
 $$
 p = \p(X\!=\!Y),\ p_j = \p(X\!=\!j),\ q_k = \sum_{j\in A}p_j^{k+1},\ q=q_2,
 $$
and let 
 $$ p_* = \max_{j\in A}p_j,\ c_+ = \log(1/q)-1,\ c_*=\log(1/p_*), $$
where $\,\log\,$ is to the base $\,1/p.$ Note that
 \beq                                                   \label{7b26}
 p_*^2<p \,,\ p^2\le q\le p_* \,.
 \eeq 
		Taking into account H\"older's inequality, we conclude that
 \beq                                                   \label{7b27}
 1\!\ge\!c_+\!\ge\!c_*\!>\!1/2. 
 \eeq
Note that $\,c_+=c_*=1\,$ if $\,\L(X)\,$ is uniform over a finite alphabet.

  Let $\,\pi_{m,n}\,$ denote a Poisson random variable with parameter $\,\l_{k,m,n}$. 
	
	 The following theorem shows that the distribution of the number of long match patterns can be well approximated by the Poisson law.

    \begin{theorem} \label{LMP} 
If $\,n\!\ge\!k\,$ and $\,m\!\ge\!k\!\ge\!1,$ then
 \beq                                                       \label{7b32}
 \d(W_{m,n};\pi_{m,n}) \le \frac{1\!-\!e^{-\l}}{\l} m'n'(2k\!+\!1)
 \(2kq_{2k} + (m'\!+\!n'\!-\!1)(p^{2k}\!+\!q^k)\). 
 \eeq \end{theorem} 

Theorem \ref{LMP} has been derived using Theorem \ref{T37} and Lemma \ref{LMP-L1}.  
		
		Denote $$ \Delta_{m,n}(k) = |\p(M_{m,n}^*\!<\!k) - \exp\(-\l\)|. $$ 

    \begin{corollary} \label{LMP-C1} 
For any constant $\,C\in\R$, as $\,m\to\infty,$ $\,n\to\infty$, 
 \beq                                                   \label{7b20}
 \max_{k\ge C+\log mn} \Delta_{m,n}(k) = O\((m\!+\!n)(mn)^{-c_+}
 (\ln mn) + (mn)^{1-2c_*}(\ln mn)^2\) \,.         
 \eeq
If $\,m\to\infty\,$ and $\,n\to\infty\,$ in such a way that $\,(\ln mn)/(\min\{m,n\}) \to 0,$ 
then
 \b                                                     \label{7b21}
 \max\limits_{1\le k\le m\wedge n} \! \Delta_{m,n}(k) =
 O\Big((m\!+\!n)(mn)^{-c_+}(\ln mn)^{1+c_+} + (mn)^{1-2c_*}(\ln mn)^{1+2c_*}\Big) \,. \e  \end{corollary}

   It is easy to see that the accuracy of estimate (\ref{7b20}) depends on the
relation between $\,m\,$ and $\,n$.  If $\,\L(X)\,$ is uniform over a finite alphabet and $\,(\ln mn)/(m\wedge n)\to 0$, then Corollary \ref{LMP-C1} implies that 
 \beq                                                   \label{7b19}
 \max\limits_{1\le k\le m\wedge n} \left| \p(M_{m,n}^*\!<\!k) - 
 e^{-\l} \right| = O\(n^{-1}(\ln n)^2\)
 \eeq 
    If $\,\L(X)\,$ is uniform over a finite alphabet and
 $$ c \le m/n \le 1/c $$  
for some constant $\,c\!>\!0$, then the right-hand side of (\ref{7b20}) becomes
$\,O(n^{-1}\ln n)$. We conject that the correct rate of convergence in (\ref{7b19}) for the uniform $\,\L(X)\,$ is $\,O\(n^{-1}\ln n\)$. 

   The reason why (\ref{7b32}) does not yield such a rate is the lack of factor
$\,e^{-\l}\,$ on the right-hand side. Results obtained for LLHR by the method of
recurrent inequalities do produce such a factor (cf. Theorem 3.12 in \cite{N11}). 

   In a more general situation one can consider NLMP with say $r$ mismatches allowed. 
   An estimate of the accuracy of Poisson approximation to the distribution of the number of long $r$-interrupted match patterns among $\,X_1,...,X_m, Y_1,...,Y_n\,$ (match patterns of length $\ge k$ with $\le r$ ``interruptions'') can be found in \cite{Mi08,N11}. 
	Neuhauser \cite{Ne94} considers a situation where $\,\L(X)\,$ may differ from $\,\L(Y)\,$ and only insertions and deletions (but no mismatches) are allowed to occur; she presents a logarithmic estimate of the rate of Poisson approximation to the distribution of the number of such long patterns.\\

    {\bf The Zubkov--Mihailov statistic}.\index{Zubkov--Mihailov statistic}
Let now $\,Y_i = X_i$ $\,(\forall i)$, $\,m=n$. Denote
 $$ N_n^* \equiv N_n^*(k) = \sum_{(i,j)\in A(n,k)} T_{ij}^*\,, $$
where $$\,A(n,k) = \{(i,j): 0\!\le\!i\!<\!j\le\!n\!-\!k\}\ \ \ (n\!>\!k).$$ 
   $\,N_n^*\,$ is the number of long match patterns in one and the same sequence, 
$\,X_1,...,X_n$. 

   Statistic $\,N^*_n\,$ was introduced by Zubkov\index{Zubkov} \&
Mihailov\index{Mihailov} \cite{ZM} who have shown that $\,\L(N_n^*)\,$ is
asymptotically Poisson $\,{\bf \Pi}(\mu)\,$ if
 $$ n^2p^k(1\!-\!p)/2 \to \mu\!>\!0,\ \ nk^tp_*^k \to 0\quad (\forall\,t\!>\!0). $$

   Note that 
 $$ M_n^* = \max\{k\!\le\!n\!: \max_{(i,j)\in A(n,k)} T_{ij}=1\} $$ 
is the {\it length of the longest match pattern}\index{LLMP} among 
$\,X_1,...,X_n$. Obviously, $$ \{M_n^*\!<\!k\}=\{N_n^*=0\}. $$ 

The next theorem evaluates the accuracy of Poisson approximation to $\,\L(N^*_n)$.

    \begin{theorem} \label{7bT14} 
If $\,n\!>\!3k\!\ge\! 3,$ then 
 $$ 
 \d(N^*_n;\pi_{n,k}^*) \le \frac{1\!-\!e^{-\lambda^*}}{\lambda^*}
 \( (n^*)^3(2k\!+\!1)\, (p^{2k}\!+\!q^k)+2(kn^*)^2q_{2k} \) + 2kn^*p^k,
 $$
where 
$\,\l^* \equiv \l_{n,k}^* = (n\!-\!3k\!+\!1)\, p^k(1+(n\!-\!3k)(1\!-\!p)/2),$ $\,n^*=n\!-\!k$, $\,\L(\pi_{n,k}^*) = {\bf\Pi}(\l^*)$.    \end{theorem} 

   Theorem \ref{7bT14} has been derived using Theorem \ref{T37} and Lemma \ref{LMP-L1}. 

    Denote $$ \Delta^*(n,k) = |\p(M_n^*<k) - \exp(-\lambda_{n,k}^*)|. $$

    \begin{corollary} \label{LMP-C2} 
As $\,n\to\infty,$ 
 \b                                               \label{7b22} 
 \max_{k\ge C+2\log n}\Delta^*(n,k) &=& O\(n^{1-2c_+}\ln n + n^{2-4c_*}
 (\ln n)^2\),\\ \label{7b23}
 \max_{1\le k<n/3}\Delta^*(n,k) &=& O\(n^{1-2c_+}(\ln n)^{1+c_+}
+ n^{2-4c_*}(\ln n)^{1+2c_*}\) . \e \end{corollary}

  If $\,\L(X)\,$ is uniform over a finite alphabet, then the right-hand side of (\ref{7b22}) is $\,O(n^{-1}\ln n),$ the right-hand side of (\ref{7b23}) is $\,O(n^{-1}(\ln n)^2)$.\\ 

    The key result behind Theorems \ref{LMP} and \ref{7bT14} is the following 

    \begin{lemma} \label{LMP-L1} 
For all natural $\,i,j,i',j'\,$ such that $\,(i,j) \neq (i',j')\,,$
 \beq                                                   \label{7b24}
 \p(T_{ij}^* = T_{i'j'}^* = 1) \le q_{2k}\,. 
 \eeq \end{lemma} 

    $\,$ 

    Denote by $$ \tau_k = \min\{n\!: N_n^*(k)\ne0\} $$ the first instance a match pattern of length $\,k\,$ appears in the sequence $\,\{X_i,i\!\ge\!1\}\,.$ Then 
  $$ \{\tau_k>n\}=\{M_n^*<k\}. $$ 
  The results on the asymptotics of $\,\tau_k\,$ can be derived from the corresponding results on $\,M_n^*\,$.\\  

  NLMP with a small number of mismatches has been considered by several authors (see \cite{Mi08,N11} and references therein). 
	
	A number of authors evaluated the accuracy of compound Poisson approximations to the distribution of NLMP (see \cite{Mi08,N11,Sch00} and references therein).\\

  \noindent{\it Open problems}.\\ 
  \noindent{3.2}. Derive uniform in $\,k\,$ estimates of (possibly shifted) Poisson approximation to $\,\L(W_{m,n})\,$ and $\,\L(N^*_n)$.\\ 
	\noindent{3.3}. Find the 2$^{\rm nd}$-order asymptotic expansions
for $\,\p(W_{m,n}\!\in\!\cdot)\,$ and $\,\p(N_n^*\!\in\!\cdot).$\\ 
  \noindent{3.4}. Check if the correct rate of convergence in (\ref{7b21}) and (\ref{7b23}) in the case of uniform $\,\L(X)\,$ is $\,O\(n^{-1}\ln n\)$.\\ 
 \noindent{3.5}. Improve the estimate of the rate of convergence in the limit theorem for the length of the longest $r$-interrupted match pattern.\\ 


							\section{Compound Poisson approximation}\label{5}

	The topic of compound Poisson (CP) approximation is vast. 
  From a theoretical point of view, the interest to the topic arises in connection with Kolmogorov's problem concerning the accuracy of approximation of the distribution of a sum of independent r.v.s by infinitely divisible laws (see \cite{AZ,LeCam65,Pre83,Pro53} and references therein). 
	
	Recall that the class of infinitely divisible distributions coincides with the class of weak limits of compound Poisson distributions \cite{Khi}. 

	The topic has applications in extreme value theory, insurance, reliability theory, patterns matching, etc. (cf. \cite{BalK,BHJ,BC01,LLR,N11}). For instance, in (re)insurance applications the sum $\,S_n = \sum_{i=1}^n\! Y_i \1\{Y_i\!>\!x_i\}\,$ of integer-valued r.v.s allows to account for the total loss from the claims $\,\{Y_i\}\,$ that exceed excesses $\,\{x_i\}$. If the probabilities $\,\p(Y_i\!>\!x_i)\,$ are small, $\,\L(S_n)\,$ can be accurately approximated by a Poisson or a compound Poisson law. 
	
   In extreme value theory one deals with the number of extreme (rare) events represented by a sum of 0-$\!$1 r.v.s (indicators of rare events). The indicators can be dependent. A well-known approach consists of grouping observations into blocks which can be considered almost independent \cite{B26}. The number of r.v.s in a block is an integer-valued r.v., hence the number of rare events is a sum of almost independent integer-valued r.v.s. 
	In all such situations the block sums are non-zero with small probabilities. 
More information concerning applications can be found in \cite{BalK,BHJ,Ge79,LLR}. 

  This section concentrates on results concerning compound Poisson (CP) approximation that can be derived from the results concerning pure Poisson approximation.

						\subsection{CP limit theorem}

	{\it Compound Poisson} (CP) distribution is the distribution of a r.v. 
  $$ \sum_{i=1}^{\pi_\l}\z_i\,, $$ 
where $\,\z_0=0,$ r.v.s $\,\pi_\l,\z,\z_1,\z_2,...$ are independent, $\,\L(\z) = {\bf\Pi}(\l)$, $\,\z_i \st \z\ (i\!\ge\!1)$. 

   We denote $\,\L(\sum_{i=1}^{\pi_\l}\z_i)\,$ by 
$\,{\bf\Pi}(\l,\z) \equiv {\bf\Pi}(\l,\L(\z))$. 

   Typically $\,\z\!\ne\!0\ w.p.\,1$. The requirement $\,\z\!\ne\!0\ w.p.\,1\,$ may be omitted. Indeed, denote $\,p=\p(\z\!\ne\!0).$ Then by Khintchin's formula (\cite{Khi}, ch. 2), 
 \beq \label{Khi} \z \st \tau_p \z', 
 \eeq 
where $\,\tau_p\,$ and $\,\z'\,$ are independent r.v.s, 
$\,\L(\z') = \L(\z|\z\!\ne\!0),$ $\,\L(\tau_p) = {\bf B}(p)$. Note that 
	$$\, {\bf\Pi}(t,\tau_p\z') = {\bf\Pi}(tp,\z') \,$$ (cf. (6.26) in \cite{N11}).  

   Let $\,\{X_{n,1},...,X_{n,n}\}_{n\ge1}\,$ be a triangle array of stationary dependent 0-$\!$1 random variables, i.e., sequence $\,X_{n,1},...,X_{n,n}\,$ is stationary for each $\,n\!\in\!\N$.   Set $$ S_n = X_{n,1}+...+X_{n,n}. $$ 
   Let $\,\z_{r,n}\,$ be a r.v. with distribution (\ref{z}). 
   The following Theorem \ref{T5.1} generalises Theorem \ref{Dep} to the case of CP  approximation. 	It states that under certain assumptions weak convergence of the cluster size distribution (see (\ref{13.3}) below) is necessary and sufficient for the CP limit theorem for $\,S_n\,$. 

   In Theorem \ref{T5.1} below we will assume (\ref{13.1}) and the following condition: 
 \beq                                         \label{13.0} 
 \limsup_{n\to\infty}\,n\p(X_{n,1}\ne0) < \infty. 
 \eeq 
Note that relation (\ref{13.1}) does not imply (\ref{13.0}) --- for example, consider the case $\,X_{n,1}\equiv X.$ Denzel \& O'Brien \cite{DO} present an example of an $\,\alpha$--mixing sequence such that (\ref{13.1}) holds though (\ref{13.0}) does not.  

 \begin{theorem} \label{T5.1} Assume conditions (\ref{13.1}), (\ref{13.0}) and $\,\Delta$. If 
 \beq \label{13.3}  \z_{r,n} \Rightarrow \zeta \qquad(n\!\to\!\infty) \eeq 
for a sequence $\,\{r\!=\!r_n\}\!\in\!{\cal R}$, then 
 \beq                                                       \label{13.2}
 S_n \Rightarrow \sum_{i=0}^{\pi(\l)} \zeta_i.
 \eeq
The limit in (\ref{13.2}) does not depend on the choice of a sequence $\,\{r_n\}\!\in\!{\cal R}$. 

   If $\,S_n\,$ converges weakly to a random variable $\,Y,$ then $\,\L(Y)\,$ is compound Poisson $\,{\bf\Pi}(\l,\z),$ where $\,\l = -\ln\p(Y\!=\!0)$. 
   If $\,\l\!>\!0,$ then (\ref{13.3}) holds for some random variable $\,\z\,$ and sequence $\,\{r\!=\!r_n\}\!\in\!{\cal R}$. \end{theorem} 

   Theorem \ref{T5.1} is effectively Theorem 5.1 from \cite{N11}.

				\subsection{Accuracy of CP approximation}\label{5.1}

   Let $\,\{X_i\}\,$ be independent r.v.s that are non-zero with small probabilities (cf. \cite{LeCam65,M87,Pre85,Z83}). Set $\,S_n := X_1+...+X_n$, and denote $$\,p_i = \p(X_i\!\ne\!0)\qquad (i\!\ge\!1).$$ 
According to Khintchin's formula (\ref{Khi}),  
	$$ X_i \st \tau_i X'_i, \eqno(\ref{Khi}^*) $$ 
where $\,\tau_i\,$ and $\,X'_i\,$ are independent r.v.s, 
$\,\L(X'_i) = \L(X_i|X_i\! \ne\! 0)$, $\,\L(\tau_i) = {\bf B}(p_i)$. Hence 
$$ S_n \st \tau_1X_1'+...+\tau_nX_n'. $$ 

	 Let $\,\z_1,...,\z_n\,$ be independent compound Poisson $\,{\bf\Pi}(p_i,X'_i)\,$ random variables. Set $\,Z_n = \sum_{i=1}^n\z_i.$ 
		Note that $\,Z_n\,$ is a compound Poisson random variable: 
		$$ \L(Z_n) = {\bf\Pi}(\l,X'_\eta), $$ 
where r.v. $\,\eta\,$ is independent of $\,X_1',...,X_n',$ 
$\,\p(\eta\!=\!j) = p_j/\l\ \ (1\!\le\!j\!\le\!n)$.
	
   A simple estimate of the accuracy of CP approximation to $\,\L(S_n)\,$ follows from the property of $\,\d\,$ and (\ref{TVD}):  
  $$ 
	\d(S_n;Z_n) \le \sum_{i=1}^n \d(\tau_i;\pi_{p_i}) \le \sum_{i=1}^n p_i^2 
	$$ (see LeCam \cite{LeCam65}, Theorem 1). 

   Zaitsev \cite{Z83} has derived an estimate of the accuracy of compound Poisson approximation that can be sharper than (\ref{LeCam}$^*$) if $\,\l=p_1+...+p_n\,$ is ``large''. The following Theorem \ref{Tzai} presents Zaitsev's result.

 \begin{theorem} \label{Tzai}  
There exists an absolute constant $\,C\,$ such that 
 \beq \label{Z83}  d_K(S_n;Z_n) \le Cp^*_n\,. \eeq \end{theorem}

   Inequality (\ref{Z83}) has been generalised to the multidimensional situation by Zaitsev \cite{Z88}.\\

   We consider now the situation where $$\,X_i'\st X'\ (\forall i).$$ 
In such a situation an estimate of the accuracy of compound Poisson approximation to $\,\L(S_n)\,$ follows from the estimate of the accuracy of pure Poisson approximation to $\,\L(\tau_1\!+\!...\!+\!\tau_n)$. 
		
	Indeed, denote  
	$$ \nu_n = \tau_1+...+\tau_n,\ \ Y = \sum_{i=1}^{\pi_\l}X'_i, $$ 
where Poisson $\,{\bf\Pi}(\l)\,$ r.v. $\,\pi_\l\,$ is independent of $\,X_1',X_2',...$. Then 
  \beq \label{H} S_n \st \sum_{i=1}^{\nu_n}X'_i. \eeq  

  It is easy to check (see, e.g., Presman \cite{Pre85}) that 
	\beq \label{Pr} \d(S_n;Y) \equiv \d\Big(\sum_{i=1}^{\nu_n}X'_i;\sum_{i=1}^{\pi_\l}X'_i\Big) \le \d(\nu_n;\pi_\l). \eeq 
	
  Kolmogorov (\cite{K56}, formula (30)) has applied (\ref{Pr}) without formulating it explicitly. Presman \cite{Pre85} was probably the first to formulate (\ref{Pr}) explicitly and present its proof.

   Presman \cite{Pre85} has evaluated $\,\d(\nu_n;\pi_\l)\,$ (and hence $\,\d(S_n;Y)$) using (\ref{Pr}) and (\ref{Pre85}). Michel \cite{M87} has applied (\ref{Pr}) and the Barbour--Eagleson estimate (\ref{be}). 
	An application of (\ref{Pr}) and (\ref{RN}) yields 
	\beq \label{Ncp} 
	\d(S_n;Y) \le 3\theta/4e + 2\delta^*\ve  + 2\delta^2. \eeq 

  According to \cite{N11}, Lemma 5.4, 
\beq d_{_G}(S_n;Y) \le d_{_G}(\nu;\pi_\l)\E|X'|. \label{13.15} \eeq 
   A combination of (\ref{2b-1}) and (\ref{13.15}) entails 
	\beq \label{cp} d_{_G}(S_n;Y) \le 
	\Big(1\wedge\frac{_4}{^3}\sqrt{2/enp}\,\Big)\l_2 \E|X'|. \eeq 

   $\,$ 

   Further results on the accuracy of compound Poisson approximation can be found in \cite{BU99,CW03,CR06,R03,Z03,Xia15}.\\ 

   \noindent{\it Open problem}.\\ 
   \noindent4.1 Evaluate constant $\,C\,$ in (\ref{Z83}).

		\subsection{CP approximation to $\,{\bf B}(n,p)$}\label{5.2}

   Below we present an estimate of the accuracy of compound Poisson approximation to the Binomial law related to the topic of pure Poisson approximation. 

   Let $\,X,X_1,...\,$ be independent Bernoulli $\,{\bf B}(p)\,$ r.v.s.  
   Presman \cite{Pre83} has shown that 
	\beq \label{Presman} 
	\sup_p\d({\bf B}(n,p);F_{n,p}) = O(n^{-2/3}), 
	\eeq 
where the compound Poisson distribution $\,F_{n,p}\,$ is constructed via {\it Poisson} distributions (a similar result in terms of $\,d_K\,$ is due to Meshalkin \cite{Me60}). 
	
  We present Presman's result in Theorem \ref{Pres} below (see also \cite{AZ}, ch. 4).\\ 

   Denote by $\,\left\lceil x \right\rceil\,$ the integer number that is the nearest to $\,x\,$ from above, and let 
	$$ \gamma = \left\lceil 3np^2\!-\!2np^3 \right\rceil,\ 
     \,\beta = \gamma\!-\!3np^2\!+\!2np^3 \!\in\![0;1),\ \,q=1\!-\!p. 
  $$ 
   Let $\,\eta_1,\eta_2,\eta_3\,$ be independent r.v.s with distributions 
 $$ \L(\eta_1) = {\bf\Pi}(pq^2\!-\!\beta/n),\ \L(\eta_2) = 
    {\bf\Pi}(p^2q\!+\!\beta/3n),\  \L(\eta_3) = {\bf\Pi}(\beta/6n). 
 $$ 
Set  $$ Y := \gamma/n\!+\!\eta_1\!-\!\eta_2 \!+\!2\eta_3. $$ 
   Note that $\,Y\,$ is a CP r.v..  One can check that 
 $$ \E Y = p,\ \E(Y\!-\!p)^2 = pq,\ \E(Y\!-\!p)^3 = pq(q\!-\!p). $$ 

  Let $\,F_{n,p}:=\L(Y_1+...+Y_n),$ where $\,\{Y_i\}\,$ are independent copies of $\,Y$.

	\begin{theorem} \label{Pres} 
There exists an absolute constant $\,C\,$ such that 
 \beq \label{Pre} 
 \d({\bf B}(n,p);F_{n,p}) \le C\ve_{n,p} \qquad(0\!\le\!p\!\le\!1/2),
 \eeq 
where $\,\ve_{n,p} =\min\big\{ np^2; p; \max\{1/(np)^{2};1/n\} \big\}.$ 
 \end{theorem}

   Bound (\ref{Presman}) follows after noticing that 
$\,\sup_{0\le p\le1/2} \ve_{n,p} = O(n^{-2/3}).$\\

   {\bf Dependent 0-$\!$1 r.v.s}. 
	Let $\,X,X_1,...\,$ be a stationary sequence of 0-$\!$1  r.v.s. The following Theorem \ref{T5.2} is an application of (\ref{Pr}) in the case of dependent r.v.s. 

   Let $\,\pi,\z_1^{(r)},\z_2^{(r)},\ldots$ be independent random variables, where $\,1 \!\le\!r\!\le\!n$, $\,\pi_{n,r}\,$ is a Poisson $\,{\bf\Pi}(kq)\,$ r.v., $\,\z_0^{(r)}=0,$ 
 \bb
 \L(\z_i^{(r)}) \!&=&\! \L(S_r|S_r\!>\!0) \quad (i\!\ge\!1),\\ 
 q \!&=&\! \p(S_r\!\ne\!0),\ \ k=[n/r]. 
 \ee  
    Denote $\,p=\p(X\!=\!1),$ 
 $$ Y_n = \sum_{i=0}^{\pi_{n,r}} \zeta_i^{(r)}. $$ 
   The distribution of $\,S_n = X_1+...+X_n\,$ can be approximated by a CP distribution $\,\L(Y_n)$.

   \begin{theorem}  \label{T5.2} \index{compound Poisson!approximation}
If $\,n\!>\!r\!>\!l\!\ge\!0,$ then
 \b                                                         \label{71}
 \d( S_n;Y_n) \!&\le&\! \kappa_{n,r}rp + (2kl+r')p + nr^{-1}\gamma_n(l),\\ 
 d_{_{G}}(S_n;Y_n) \!&\le&\! rp\min\!\Big\{np\,; \frac{_4}{^3}\sqrt{2np/e}\,\Big\}
 + (2kl+r')p + n\gamma_n(l),                        \label{G}
 \e
where $\,r'\!=\!n\!-\!rk,$ $\,\kappa_{n,r} = \min\{ 1\!-\!e^{-np}\,; 3/4e+(1\!-\!e^{-np})rp \}\,$ and $\,\gamma_n(l) = \min\{4\a(l)\sqrt{r}\,;\beta_n(l)\}$. \end{theorem}

  Theorem \ref{T5.2} is effectively Theorem 5.2 from \cite{N11}. 
	
	If the random variables $\,\{X_i\}\,$ are independent, then
(\ref{71}) with $\,r\!=\!1,$ $\,l\!=\!0\,$ yields (\ref{be}) and (\ref{RN}). 

   If the random variables $\,\{X_i\}\,$ are $m$--dependent, then one can choose $\,l\!=\!m$, $\,r=\lceil\sqrt{mn}\,\rceil$, the smallest integer greater than or equal to $\,\sqrt{mn}\,,$ and get the estimate $\,\d( S_n;Y_n) \le 4p\lceil\sqrt{mn}\,\rceil.$ 
	
	Further reading on the topic of the accuracy of compound Poisson approximation to the distribution of a sum of dependent r.v.s includes \cite{Roos94,ChVe10} and references therein.\\ 

 \noindent{\it Open problem}.\\ 
 \noindent4.2. Evaluate constant $\,C\,$ in (\ref{Pre}).\\ 
	

					\section{Poisson process approximation}\label{4}

  The topic of point process approximation is vast; an interested reader is referred to \cite{DV08,MKM}. 
  This section concentrates on the results concerning Poisson process approximation that are closely related to the results on Poisson approximation  to the distribution of a sum of 0-$\!$1 random variables.\\ 

  \indent{\bf Point process counting locations of rare events}. 
   Let $\,\{\xi_i,i\!\ge\!1\}\,$ be Bernoulli r.v.s (e.g., $\,\xi_i = \1\{X_i\!>\!u_n\}$, where $\,u_n\,$ is a ``high'' level). Then 
 \beq                              \label{2b1} 
 S_n(\cdot) = \sum_{i=1}^n \1\{\,i/n\!\in\!\cdot\} \xi_i 
 \eeq 
can be called a ``Bernoulli process''\index{Bernoulli process}. 

   $S_n(\cdot)\,$ counts {\it locations} of extreme/rare events represented by r.v.s $\,\{\xi_i\}$. 
	A typical example of a rare event is an exceedance of a high threshold. 

  For instance, let $\,X,X_1,X_2,...$ be a stationary sequence of random variables, and let $\,\{u_n\}\,$ be a sequence of levels. Set $\,\xi_i = \1\{X_i\!>\!u_n\}$. Then $\,S_n(\cdot) = N_n(\cdot,u_n),$ where 
 $$  N_n(B,u_n) = 
 \sum_{i=1}^n \1\{\,i/n\!\in\!B, X_i\!>\!u_n\} \qquad(B\!\subset\!(0;1]).
 \eqno(\ref{2b1}^*) $$ 
   Process $\,N_n(\cdot,u_n)\,$ counts {\it locations} of exceedances of level $\,u_n$. 

  Let $\,\{r\!=\!r_n\}\,$ be a sequence obeying (\ref{rr}). 
We denote by $\,\z_{r,n}\,$ a r.v. with distribution (\ref{z}).

  \begin{theorem} \label{2bt1} 
Assume (\ref{13.1}), (\ref{13.0}) and mixing condition $\,\Delta$. 
   If (\ref{1}) holds, then 
 \beq                                                      \label{2b0}
 N_n(\cdot,u_n) \Rightarrow N(\cdot), \eeq 
where $\,N(\cdot)\,$ is a Poisson point process with intensity rate $\,\l$. 
  \end{theorem}

	  Theorem \ref{2bt1} is a particular case of Theorem 7.2 in \cite{N11}. The necessity part of Theorem \ref{2bt1} is given by Theorem \ref{Dep}: if (\ref{2b0}) holds, then so does (\ref{1}).  
	Leadbetter et al. \cite{LLR}, Theorem 5.2.1, present a version of Theorem \ref{2bt1} with condition $\,(D')\,$ instead of (\ref{1}).\\ 
	
	Denote by $\,\Xi_n\,$ a Poisson point process with intensity measure 
	$$\,\l(\cdot) = \sum_{i=1}^n p_i\1\{i/n\!\in\!\cdot\},$$ 
where $\,p_i=\p(\xi_i\!=\!1)$. 

   The accuracy of Poisson process approximation to $\,\L(S_n(\cdot))\,$ has been evaluated by Brown \cite{Br83} and Kabanov et al. \cite{KLSh}, Theorem 3.2: if $\,\{\xi_i\}\,$ are independent, then 
  $$ \d(S_n(\cdot);\Xi_n(\cdot)) \le \sum_{i=1}^n p_i^2. \eqno(\ref{LeCam}'') $$ 

  Arratia et al. \cite{AGG} have generalised (\ref{LeCam}$''$) to the case of dependent Bernoulli r.v.s. 
	
	Ruzankin \cite{R04} and Xia \cite{Xia05} present  estimates of the accuracy of Poisson process approximation in terms of a $\,d_{_{G}}$-type distance. 
	
   In the general case (when the limiting distribution of $\,\z_{r,n}\,$ is not degenerate) the limiting distribution of $\,N_n(\cdot,u_n)\,$ is necessarily compound Poisson (Hsing et al. \cite{HHL}, see also \cite{N11}, ch. 7).\\

					{\bf Excess process}.  
	 Let $\,X,X_1,X_2,...\,$ be a stationary sequence of r.v.s. 

  If one is interested in the joint distribution of exceedances of several levels among $\,X_1,...,X_n$, a natural tool is the excess process $\,N_n^\ve(\cdot)$. 	We give the definition of the excess process below. 
	
   Suppose there is a sequence $\,\{u_n(\cdot),n\!\ge\!1\}\,$ of functions on $\,[0;\infty)\,$ such that function $\,u_n(\cdot)\,$ is strictly decreasing for all large enough $\,n,$ $\,u_n(0)=\infty,$  
 \b                                                    \label{2b9}
 && \limsup_{n\to\infty} n\p(X\!>\!u_n(t))<\infty \qquad\qquad 
 (0\!<\!t\!<\!\infty),\\                            		\label{2b5} 
 && \lim_{n\to\infty} \p(M_n\!\le\!u_n(t)) = e^{-t}\quad\qquad\quad\ \ \ (t\!\ge\!0),
 \e 
where $\,M_n = \max\{X_1,...,X_n\}\,$ is the sample maximum. 
   Conditions (\ref{2b9}) and (\ref{2b5}) mean that $\,u_n(\cdot)\,$ is a ``proper'' normalising sequence\index{normalising sequence} for the sample maximum. 

   Set $$\,N_n^\ve(t) = \sum_{i=1}^n \1\{X_i\!>\!u_n(t)\},$$ 
where $\,t\!>\!0.$ Given $\,B\!\subset\![0;\infty),$ we call $\,\{N_n^\ve(t),t\!\in\!B\}\,$ the {\it excess process}. 

  Process $\,N_n^\ve(\cdot)\,$ describes variability in the {\it heights} of the extremes. 
	
	Note that $\,N_n^\ve(\cdot)\,$ is the ``tail empirical process'' for $\,Y_{n,1},...,Y_{n,n},$ where $\,Y_{n,i}\!=\!u_n^{-1}(X_i)$: 
	\beq \label{2be} N_n^\ve(t) = \sum_{i=1}^n \1\{Y_{n,i}\!<\!t\}. \eeq 

  There is considerable amount of research on the topic of tail empirical processes (see, e.g., \cite{DIK} and references therein). 
	
	Below we present necessary and sufficient conditions for the weak convergence of the excess process to a Poisson process in Theorem \ref{2bt3} (cf. \cite{N11}, ch. 7).\\ 

 
   First, we recall the definitions of mixing (weak dependence) conditions. 

   Given $\,0\!<\!t_1\!<\!...\!<\!t_k\!<\!\infty,$ where $\,k\!\ge\!1,$ and a sequence $\,\{u_n(\cdot)\}_{n\ge1},$ we denote $$\,\t=(t_1,...,t_k),\ \,u_n(\t)=(u_n(t_1),...,u_n(t_k)).$$ 
   Let $\,{\cal F}_{l,m}(\t)\,$ be the $\,\sigma$--field generated by the events 
$\,\{ X_i\!>\!u_n(t_j)\},$ $\,l\!\le\!i\!\le\!m, 1\!\le\!j\!\le\!k$; mixing (weak dependence) coefficient $\,\a_n(l_n) := \a\(l_n,u_n(\t)\)\,$ is defined as above.\\ 
  
 
   {\it Condition $\,\Delta(\{u_n(\t)\})\,$} is said to hold if 
$\,\a_n(l_n) \to 0 \,$ for some sequence $\,\{l_n\}\,$ 
such that $\,l_n\to\infty,\ l_n/n\to 0\,$ as $\,n\to\infty$.\\ 

    {\it Condition $\,\Delta^*\,$} holds if $\,\Delta(\{u_n(\t)\})\,$ is in force 
$\,(\forall\, 0\!<\!t_1\!<\!\ldots\!<\!t_k\!<\!\infty,$ $\,k\!\ge\!1).$\\ 

   {\it Class} $\,{\cal R}(\t)$. If $\,\Delta\{u_n(\t)\}\,$ holds, then there exists a sequence $\,\{r_n\}\,$ such that (\ref{rr}) holds 
(for instance, one can take $\,r_n = \big[\sqrt{n\max\{l;n\a_n(l_n)\}}\,\big]$). 
We denote by $\,{\cal R}(\t)\,$ the class \label{classR} of all such sequences.\\  

  The next condition describes the joint distribution of exceedances of several levels.\\ 

  We say that {\it condition $\,C_{\t}'\,$} holds if there exists a sequence 
$\,\{r_n\}\!\in\!{\cal R}(\t)\,$ such that for every $\,1\!\le\!i\!<\!j\le k\,$ and every $\,t_i\!<\!t_j\,$ from $\{t_1,...,t_k\}\,$\\ 
   \noindent(a) 
 $\, \p(N_r[u_n(t_{i-1}); u_n(t_i))\!=\!1) \sim\! \frac rn (t_i\!-\!t_{i-1}),$ 
 $\, \p(N_r[u_n(t_{i-1}); u_n(t_i))\!=\!j) =  o(\frac rn)\ \ (j\!\ge\!2),$ 

   \noindent(b)  
  $\, \p(N_r(u_n(t_i))\!>\!0, N_r[u_n(t_i);u_n(t_j))\!>\!0) = o(r/n).$\\ 

   {\it Condition $\,C'\,$\index{condition $\,C'$}} holds if {\it $\,C_{\t}'\,$ is valid for all} $\,0\!<\!t_1\!<\!...\!<\!t_k\!<\!\infty,$ $\,k\!\ge\!1$.\\

    \begin{theorem} \label{2bt3}              
Assume mixing condition condition $\,\Delta,$ and let $\,\pi(\cdot)\,$ denote a Poisson process with intensity rate 1. Then
 \beq                                                   \label{2b7}
 N_n^\ve(\cdot) \Rightarrow \pi(\cdot) \eeq 
if and only if condition $\,C'\,$ holds. \end{theorem} 

 $\,$ 

   \noindent{\bf Example 5.1}.  
	Let $\,X,X_1,X_2,...\,$ be i.i.d.r.v.s with the distribution function (d.f.) $\,F$. 
	Denote $\,K^*=\sup\{x\!: F(x)\!<\!1\},$ and assume that 
 $$ 
 \p(X\!\ge\!x)/\p(X\!>\!x) \to 1 \eqno(G) 
 $$ 
as $\,x\to K^*\,$ (Gnedenko's condition\index{Gnedenko's condition} \cite{G43}). 
	Set $\,u_n(t) = F_c^{-1}(t/n)$, where $$\,F_c(\cdot):=\p(X\!>\!\cdot).$$ 
Then excess process $\,\{N_n^\ve(\cdot),t\!\in\![0;1]\}\,$ converges weakly to a pure Poisson process $\,N\,$ with intensity rate 1. Process $\,N\,$ admits the representation  $$ N \st \sum_{j=1}^{\pi(1)} \gamma_j(\cdot), $$ 
where $\,\gamma_j(t) \st \1\{\xi\!<\!t\}\,$ and r.v. $\,\xi\,$ has uniform $\,{\bf U}[0;1]\,$ distribution. \hspace*{\fill} $\Box$\\

  The accuracy of approximation $\,N_n^\ve(\cdot) \approx N(\cdot)\,$ can be evaluated as well (cf. Deheuvels \& Pfeifer \cite{DP88}, Kabanov \& Liptser \cite{KL83}, Novak \cite{N11}, ch. 8). 
	
		Note that (\ref{Pr}) is applicable. Given $\,T\!>\!0,$ let $\,\pi(np)\,$ denote a Poisson $\,{\bf\Pi}(np)\,$ r.v., where $\,p=\p(X\!>\!u_n(T))$. 
  Let $\,\eta,\eta_1,\eta_2,...$ be independent of $\,\pi(np)\,$ i.i.d. processes with the distribution 
	$$ \L(\eta(\cdot)) = \L(\1\{X\!>\!u_n(\cdot)\}|X\!>\!u_n(T)) 
	\equiv \L(\1\{Y_{n,1}\!<\!\cdot\}|Y_{n,1}\!<\!T) 
	$$ 
($i\!\ge\!1$). An application of (\ref{Pr}) and (\ref{RN}) yields 
	\beq \label{est}  
	\d\Big( N_n^\ve(\cdot);\sum_{i=1}^{\pi(np)} \eta_i(\cdot) \Big) 
     \le 3p/4e + 2(1\!-\!e^{-np})p^2\ve \!+\! 2(1\!-\!e^{-np})^2p^2, 
	\eeq  
where $\,\ve = \min\!\big\{1; \(2\pi[(n\!-\!1)p]\)^{-1/2} + 2(1\!-\!e^{-np})p/(1\!-\!1/n)\big\}\,$ (\cite{N11}, Theorem 8.3). 

  Note that $\,\sum_{i=1}^{\pi(np)} \eta_i(\cdot)\,$ is a Poisson process. 
	If $\,F_c\,$ is a continuous decreasing function, then $\,\eta(\cdot)\st\1\{\xi\!<\!\cdot\},$ where $\,\L(\xi) = {\bf U}[0;1].$\\ 
	
	Let $\,X,X_1,X_2,...\,$ be i.i.d.r.v.s, and let $\,B\!\subset\![0;\infty)\,$ be a closed set. According to (6.5) in \cite{DP88} and (\ref{tdvM}), the total variation distance between $\,\{\sum_{i=1}^n \1\{Y_{n,i}\!<\!t\},\,t\!\in\!B\}\,$ and the approximating Poisson process coincides with $\,\d({\bf B}(p);{\bf\Pi}(p)),$ where $\,p=\p(Y_{n,1}\!\in\!B)$. 

  In a general situation excess process $\,\{N_n^\ve(\cdot)\}\,$ may converge weakly to a process of more complex structure: 
 \beq \label{2b40}    \{N_n^\ve(t),\,t\!\le\!T\} \Rightarrow 
 \Bigg\{ \sum_{j=1}^{\pi(T)} \gamma_j(t/T),\, t\!\le\!T \Bigg\}, 
 \eeq  
where $\,\pi(T)\,$ is Poisson $\,{\bf\Pi}(T)$, $\,\{\gamma_j(\cdot)\}\,$ are independent jump processes. 
	
   Process $\,\Big\{\sum_{j=1}^{\pi(T)} \gamma_j(t)\Big\}\,$ 
can be called {\it Poisson cluster process}\index{Poisson!cluster process} or {\it compound Poisson process of the second order} (regarding the standard CP process as a ``compound Poisson process of the first order''). 

   Necessary and sufficient conditions for the weak convergence of the excess process to a compound Poisson process or a Poisson cluster processes are presented in \cite{N11}, ch. 7, 8. For an estimate of the total variation distance between $\,N_n^\ve(\cdot)\,$ and the approximating process in the case of weakly dependent r.v.s see \cite{N11}, Theorem 8.3.\\

	 {\bf General point process of exceedances}.  
   Consider now a two--dimensional point process $\,N^*_n\,$ that counts locations of rare events (e.g., exceedances of ``high'' thresholds) as well as their ``heights''$\!$: for any Borel set $\,A\!\subset\!(0;1]\!\times\![0;\infty)\,$ we set 
  \beq                                                           \label{2b11}
  N^*_n(A) := \sum_{i=1}^n \1\{\,\(i/n,u_n^{-1}(X_i)\) \in A \,\}. 
  \eeq 
  If $\,\{X_i\}\,$ are i.d.d.r.v.s, or if $\,\{X_i,i\!\ge\!1\}\,$ is a strictly stationary sequence obeying certain mixing conditions, then 
$\,N^*_n(\cdot)\,$ converges weakly to a pure Poisson point process (Adler\index{Adler} \cite{A78}). Theorem \ref{2bc} below presents Adler's result. 

   We will need a multilevel version of the ``declustering'' {\it condition $(D')$}\index{condition $(D')$}: 
 $$ \lim_{n\to\infty} n\sum_{i=1}^r \p(X_{i+1}\!>\!u_n(t), X_1\!>\!u_n(t)) = 0
 \eqno{(D'_+)} $$  
for any sequence $\,\{r\!=\!r_n\}\in{\cal R}(t),\, 0\!<\!t\!<\!\infty.$

   \begin{theorem}  \label{2bc} 
   If conditions $\,\Delta\,$ and {\rm(}$D'_+${\rm)} hold, then $\,N^*_n\,$ converges weakly to a pure Poisson point process $\,N^*\,$ on $\,(0;1]\!\times\![0;\infty)\,$ with the Lebesgue intensity measure. \end{theorem}

 \noindent{\bf Example 5.2.} 
Let $\,Y,Y_1,Y_2,...\,$ be a sequence of i.i.d.r.v.s with exponential {\bf E}(1) distribution, and set  $$ X_i=Y_i+Y_{i+1}\,. $$ 
Evidently, $\,\{X_i,i\!\ge\!1\}\,$ is a stationary sequence of 1--dependent r.v.s. 

   Let $\,u\equiv u_n(t)=\ln[t^{-1}n\ln n],$ $\,t\!>\!0.$ Then $\,\p(X\!>\!u_n(t)) \sim t/n,$ and condition ($D'_+$) holds. 
	According to Theorem \ref{2bc}, $\,N_n^*\Rightarrow$ $N^*,$ the Poisson point process with the Lebesgue intensity measure (cf. \cite{N11}, ch. 7). 
	\hspace*{\fill} $\Box$\\

   Adler's result has been generalised to the case of compound Poisson approximation: necessary and sufficient conditions for the weak convergence of $\,N^*_n\,$ to a compound Poisson point process can be found in \cite{N11}, ch. 7. 
	Necessary and sufficient conditions for the weak convergence of $\,N^*_n\,$ to a Poisson cluster process are given in \cite{N11}, ch.	8. 
		
	An estimate of the accuracy of approximation $\,N_n^*(\cdot) \approx \sum_{j=1}^{\pi(T)} \gamma_j(\cdot)\,$ in terms of the $\,d_{_{G}}(X;Y)$-type distance is given in \cite{BNX02}.\\ 

  \noindent{\it Open problem}.\\ 
  \noindent5.1. Improve the estimate of the accuracy of approximation $\,N^*_n \approx N^*\,$presented in \cite{BNX02}.\\ 
	
	\section*{Acknowledgements} 
	
The author is grateful to P.S.Ruzankin and the referee for helpful remarks.\\ 
   
	     \end{document}